\documentclass[a4paper, 12pt]{article}
\usepackage{mathrsfs}
\usepackage{amsmath,amsthm,amssymb,amsfonts,amscd,amstext}
\usepackage{latexsym}
\usepackage{amssymb}
\usepackage{xspace}
\usepackage[utf8]{inputenc}
\usepackage{hyperref}
\usepackage{enumerate}
\usepackage{bm}
\usepackage{bbm}
\usepackage{verbatim}
\usepackage{mathpazo}
\usepackage{xcolor}
\usepackage[sort&compress,numbers,square]{natbib}
\theoremstyle{plain}
\newtheorem{theorem}{Theorem}[section]

\newtheorem{proposition}[theorem]{Proposition}
\newtheorem{lemma}[theorem]{Lemma}
\theoremstyle{definition}
\newtheorem{definition}[theorem]{Definition}
\newtheorem{remark}[theorem]{Remark}

\numberwithin{theorem}{section}
\numberwithin{equation}{section}

\newcommand{\cA}{{\mathcal A}}

\newcommand{\cC}{{\mathcal C}}

\newcommand{\cE}{{\mathcal E}}

\newcommand{\cN}{{\mathcal N}}

\def\R{\mathbb{R}}

\newcommand{\eps}{\varepsilon}

\newcommand{\dis}{\displaystyle}
\newcommand{\todo}[1]{ {\tt [$\clubsuit$ #1 $\clubsuit$]}}
\usepackage{algorithm, algpseudocode} % use algorithm and algorithmicx for typesetting algorithms
\usepackage{mathrsfs} % for \mathscr command

\usepackage{lipsum}

\title{Bifurcation results and multiple solutions for the fractional $(p,q)$-Laplace operators}
\author{Emmanuel Wend-Benedo Zongo and Pierre Aime Feulefack}
\date{}
\begin{document}
\maketitle
\begin{abstract}
\noindent
We investigate a nonlinear nonlocal eigenvalue problem involving the sum of fractional $(p,q)$-Laplace operators $(-\Delta)_p^{s_1}+(-\Delta)_q^{s_2}$ with $s_1,s_2\in (0,1)$; $p,q\in(1,\infty)$ and  subject to  Dirichlet boundary conditions in an open  bounded set of $\mathbb{R}^N$. We prove bifurcation results from trivial solutions and from infinity for the considered nonlinear nonlocal eigenvalue problem. We also show the existence of multiple solutions of the nonlinear nonlocal problem using variational methods.\\
%\noindent\todo{To add after thing at the end}.\\
\\

\noindent
Keywords: Nonlocal operators, nonlinear eigenvalue values, bifurcation results, critical point theory, fractional $(p,q)$-Laplacian, multiple solutions.\\
\\
AMS 2010 subject classification: 35J20, 35J92, 35P30, 35R11, 49J35, 49J40.
\end{abstract}
\maketitle
\tableofcontents
\bigskip
%\todo{Avant tout propos il faut mentionner que nous deux nous souffrons. C'est dans la souffrance que cet article est entrain d'être réalisé. Les affaires ne marchent pas du tout mais on croit au miracle de Dieu notre créateur.}
  \section{Introduction and main results}
Let $\Omega\subset\mathbb{R}^N$, $N\ge 2$  be an open bounded set with Lipschitz boundary. In this article, we are interested in the study of the following nonlinear eigenvalue problem involving a sum  of mixed fractional $(p,q)$-Laplace operators subject to the Dirichlet boundary condition,
\begin{equation}\label{Eq1}
	\begin{split}
	\quad\left\{\begin{aligned}
(-\Delta)_p^{s_1}u+(-\Delta)_q^{s_2}u&= \lambda |u|^{q-2}u && \text{ in \quad $\Omega$}\\
		u &=0     && \text{ on } \quad \R^N\setminus\Omega,
	\end{aligned}\right.
	\end{split}
	\end{equation}
where,  $s_1,s_2\in (0,1)$, $\lambda\in \R$ and $p,q\in (1,+\infty)$. For $s\in (0,1)~\text{and}~r\in (1,+\infty)$, the operator $(-\Delta)^s_r$ stands for the fractional $r$ -Laplacian, and it is defined up to a normalization constant and compactly supported smooth functions $u:\R^N\to \R$, by
$$(-\Delta)^s_ru(x):=P.V~\int_{\R^N}\frac{|u(x)-u(y)|^{r-2}(u(x)-u(y))}{|x-y|^{N+rs}}\ dx,$$ 
where, $P.V$ stands for the Cauchy principle value. This covers the usual definition of the fractional Laplacian $(-\Delta)^s$ when $r=2$~ (see \cite{kwasnicki17}). In particular, for $s_1=s_2=1,$  problem (\ref{Eq1}) reduces to the following classical eigenvalue problem involving the $(p,q)$-Laplace operators,
% \begin{equation}\label{e2}
%     -\Delta_pu-\Delta_q u=\lambda |u|^{q-2}u~~\text{in}~~\Omega,~~~~~~~~~~~u=0~~~~\text{in}~~\partial\Omega,
% \end{equation}
\begin{equation}\label{e2}
	\begin{split}
	\quad\left\{\begin{aligned}
-\Delta_pu-\Delta_q u&= \lambda |u|^{q-2}u && \text{ in \quad $\Omega$}\\
		u &=0     && \text{ on } \quad \partial\Omega,
	\end{aligned}\right.
	\end{split}
	\end{equation} 
where formally, the operator $\Delta_ru:=\text{div}(|\nabla u|^{r-2}\nabla u)$ denotes the classical $r$-Laplacian. This kind of problem has been widely studied in the literature and occurs in quantum field theory, plasma physics, and chemical reaction design. Among the results related with problem \eqref{e2}, we refer the reader to  the articles \cite{Fran, Fras, F18,bobkov23} and  the references therein.  The special case $s_1=s_2=1$ and $q=2$, namely, the nonlinear eigenvalue problem
\begin{equation}\label{ZRq2}
	\begin{split}
	\quad\left\{\begin{aligned}
-\Delta_pu-\Delta u&= \lambda u && \text{ in \quad $\Omega$}\\
		u &=0     && \text{ on } \quad \partial\Omega,
	\end{aligned}\right.
	\end{split}
	\end{equation} 
has been considered in \cite{zongo22}. 
We note that problem (\ref{Eq1})  is a nonlocal counterpart of problem (\ref{e2}) and \eqref{ZRq2}  and, generalized in a  sense, problem (\ref{e2}) and \eqref{ZRq2} for $s_1,s_2\in (0,1)$ and $p,q\in (1,\infty)$. This is  due to the asymptoitics limit $(-\Delta)^s_ru\to -\Delta_ru$ as $s\to 1^-$,  for compactly supported smooth functions $u$ (see \cite{brasco15}). 

Fuel by various  concrete applications in many fields of sciences such as,  finance, phase transitions, stratified materials, anomalous diffusion, crystal dislocation, soft thin films, semipermeable membranes,  conservation laws, ultra-relativistic limits of quantum mechanics \cite{molica15}; the study of nonlocal elliptic problems driven by the sum of two nonlocal operators has gained a tremendous popularity in the last years, see \cite{ghanmi22,zhi20,biswas22,ambrosio22,farcaseanu17,nguyen20,razani23,bhakta19,goel20} and the references therein. This is due to their ability to describe the superstition of two processes with different scales, which finds wide applications in biological population in an ecological system \cite{dipierro21}, and in many other fields of sciences and engineering.

The main purpose of this  paper is  to extend the results proven in    \cite{ZB}, we will  analyse the spectrum of equation (\ref{Eq1})  (see also \cite{Bah,Lin,Fran, Fras, F18,bobkov23} for the same kind of results),  prove the existence of bifurcation branches forced by the eigenvalues of the homogeneous fractional $q$-Laplacian,
\begin{equation}\label{diric}
	\begin{split}
	\quad\left\{\begin{aligned}
(-\Delta)^{s_2}_q u&= \lambda  |u|^{q-2}u && \text{ in \quad $\Omega$}\\
		u &=0     && \text{ on } \quad \R^N\setminus\Omega,
	\end{aligned}\right.
	\end{split}
	\end{equation}
and prove  the  existence of multiple solutions of  problem \eqref{Eq1} using variational methods and some  topological techniques.  

 We will rely on the methods used in \cite{ZB} for problem driven by the sum of classical $(p,q)$-Laplace  operators to analyse the  spectrum and study the existence bifurcation branches of problem (\ref{Eq1}). Moreover, for $s_1,s_2\in (0,1)$ and $p,q\in (1,\infty)$, nonlinear eigenvalue problem involving the sum of two nonlocal operators $(-\Delta)^{s_1}_p+(-\Delta)^{s_2}_q$ as in equation \eqref{Eq1} has  been  studied in the literature and   full description of the set of eigenvalues  is provided, see  \cite{farcaseanu17,nguyen20} and the references therein, and \cite{biswas22} for more general eigenvalue problems of fractional $(p,q)$-Laplace operators with two parameters. 
 
 We would like to point out that in the last several years, many authors have dedicated special attention  on the study of (global) bifurcation type results for  elliptic problems involving the classical $p$-Laplacian  subject to the Dirichlet boundary condition \cite{AM,rabinowitz71},
 \begin{equation}\label{Claasical-p}
	\begin{split}
	\quad\left\{\begin{aligned}
-\Delta_p u&= f(\lambda, x, u) && \text{ in \quad $\Omega$}\\
		u &=0     && \text{ on } \quad \partial\Omega,
	\end{aligned}\right.
	\end{split}
	\end{equation} 
where the nonlinear perturbation $f: \R\times\Omega\times\R\mapsto\R$ satisfies the Carath\'eodory condition for the  second and third variables with  $f(\lambda,x,0)=0$ and $\Omega$ is a bounded domain with Lipschitz boundary. For the motivations that lead to the novelity of the bifurcation results  of problem \eqref{Eq1} studied in this paper, we point out that  analogue bifurcation results for problem \eqref{Claasical-p}  have been   successfully extended to the nonlocal counterpart in the last ten years by many authors \cite{MCP,del16,wang23},  
 \begin{equation}\label{Fractionl-p}
	\begin{split}
	\quad\left\{\begin{aligned}
(-\Delta)_p^{s} u&= f(\lambda, x, u) && \text{ in \quad $\Omega$}\\
		u &=0     && \text{ on } \quad \R^N\setminus\Omega,
	\end{aligned}\right.
	\end{split}
	\end{equation}
 where the Dirichlet condition is imposed on the complement of $\Omega.$ In the particular case $p=2$, we refer to paper  \cite{MCP}, where the authors used bifurcation theory to study  the existence of continua, bifurcating from infinity at $\lambda_1$, the first eigenvalue of $(-\Delta)^s$, containing positive and negative solutions. They also provided sufficient conditions for each of the continuum to bifurcate to the left or right of the hyperplane $\lambda=\lambda_1$ in $\R\times L^{\infty}(\Omega)$ and they also discussed multiplicity of solutions near the principal eigenvalue $\lambda_1.$
 
 However, to our best knowledge, there is no literature related to the bifurcation discussion for problem like  \eqref{Eq1},
involving  the sum of two nonlocal operators $(-\Delta)^{s_1}_p+(-\Delta)^{s_2}_q$. Following the local case approach as in \cite{ZB} , we provide a description of the set of eigenvalues of  problem \eqref{Eq1},  prove the existence of bifurcations branches for equation (\ref{Eq1}) emanating from the eigenvalues of the $q$-Laplacian  and prove multiplicity results. \\
%We not that  every weak nonnegative solution of \eqref{Eq1} is bounded \cite{biswas22}.\\

% In \cite{MCP}, the authors used bifurcation theory to study the following nonlocal problem:
% \begin{equation}\label{e1}
% 	\begin{split}
% 	\quad\left\{\begin{aligned}
% (-\Delta)^{s} u&= \lambda u+f(\lambda, x, u) && \text{ in \quad $\Omega$}\\
% 		u &=0     && \text{ on } \quad \R^N\setminus\Omega,
% 	\end{aligned}\right.
% 	\end{split}
% 	\end{equation}
% % \begin{equation}\label{e1}
% % \begin{cases}
% % (-\Delta)^{s} u=\lambda u+f(\lambda, x, u) ~~~\text{in}~\Omega,\\
% % u=0~~~~~~~~~~~~\text{in}~\R^N\backslash\Omega,
% % \end{cases}
% % \end{equation}
% where $$(-\Delta)^su(x)= P.V.~\int_{\R^N}\frac{u(x)-u(y)}{|x-y|^{N+2s}}\ dy$$ for fixed $s\in (0,1),$ $\Omega\subset \R^N~(N>2s)$ is a bounded domain with $C^{1,1}$ boundary and $\lambda\in\R$ is the bifurcation parameter. The nonlinear perturbation $f: \R\times\Omega\times\R\mapsto\R$ is the Carathéodory function. In \cite{MCP}, they  
% In this paper, we discuss the asymptotic behaviour of the the spectrum of the nonlinear eigenvalues problem \eqref{Eq1}.\\

%\indent For $q \neq 2$, we prove for equation \eqref{Eq1}, the existence of bifurcation branches (forced by the fractional $p$-Laplace operator) from the eigenvalues of a {\it nonlinear, but $q$-homogenous} equation. 

%\textbf{Il faut qu'on mentionne l'article 11 qui est un article phare dans l'étude de la bifurcation du p-fractional, dont la preuve est basée sur le calcul du dégré topologique de Leray-Schauder}

Our  first  result writes as follows:
\begin{theorem} \label{MainResult}
It holds:

\begin{itemize}
\item[1)] For every fixed $\rho > 0$ there exists a sequence of eigenvalues $\big(\lambda_k(s_1, s_2, p,q;\rho)\big)_k$ with corresponding eigenfunctions $\pm u_k$ satisfying  $\int_\Omega |u_k|^qdx = \rho$, with 
$$\text{$\lambda_k(s_1, s_2, p,q;\rho) \to +\infty$ \quad as\quad $k \to \infty$.}$$
%\par \medskip \noindent
\item[2)] The variational eigenvalues $\lambda_k(s_2,q)$ of problem \eqref{diric} are bifurcation points from $0$ if $0<s_2<s_1<1<q<p<\infty$, and, are  bifurcation points from infinity if  $0<s_1<s_2<1<p<q<\infty$, for the nonlinear eigenvalues $\lambda_k(s_1, s_2, p,q;\rho)$.
%\par \medskip \noindent
\item[3)] For any fixed $\lambda \in (\lambda_k(s_2,q),\lambda_{k+1}(s_2, q))$ there exist $k$ eigenvalues of \eqref{Eq1} with  $$\lambda = \lambda_1(s_1, s_2, p,q;\rho_1)= \dots = \lambda_k(s_1, s_2, p,q;\rho_k),$$ 
with corresponding eigenfunctions $\pm u_k$, satisfying $\int_\Omega |u_k|^q = \rho_k.$
\end{itemize}
\end{theorem}

Our second result on multiple solutions writes as follows:

\begin{theorem}\label{Mainresult2}
Assume that $0<s_2<s_1<1<q<p<\infty$ or $0<s_1<s_2<1<p<q<\infty$, and suppose that $\lambda\in(\lambda_k(s_2,q),\lambda_{k+1}(s_2,q)),$ $k\ge 1$. Then equation (\ref{Eq1}) has at least $k$ pairs of nontrivial solutions.
\end{theorem}

Note that by the symmetry of the parameters $s_1,s_2\in (0,1)$ and $p,q\in (1,\infty)$, our techniques still apply if one interchanges the role of the couple $(s_1,p)$ and $(s_2,q)$. More precisely, the results of Theorem \ref{MainResult} and Theorem \ref{Mainresult2} hold for the problem:
  \begin{equation}\label{Eq2-p}
	\begin{split}
	\quad\left\{\begin{aligned}
(-\Delta)_q^{s_2} u+(-\Delta)_p^{s_1}u&= \lambda |u|^{p-2}u && \text{ in \quad $\Omega$}\\
		u &=0     && \text{ on } \quad \R^N\setminus\Omega.
	\end{aligned}\right.
	\end{split}
	\end{equation}

We comment on the proof of Theorem \ref{MainResult}. 
The proof of the first item of the Theorem \ref{MainResult} is presented in Section \ref{Spec}, where we analyse the set of eigenvalues and the existence of eignfunctions of problem \eqref{Eq1}, depending whether, $0<s_2<s_1<1<q<p<\infty$ or  $0<s_1<s_2<1<p<q<\infty$. In the case $0<s_2<s_1<1<q<p<\infty$, the  functional associated to problem \eqref{Eq1} is coercive and we use the  Direct Method in the Calculus of Variations in order to find critical points of the associated energy functional. In the case $0<s_1<s_2<1<p<q<\infty$, the energy functional associated to problem \eqref{Eq1} is not  coercive. To prove that the associated energy functional has a critical point in $W^{s_2,q}_0 (\Omega) \backslash \{0\}$, we constrain the functional on the Nehari manifold and show through a series of propositions that the critical point in the Nehari manifold is in fact a solution of problem \eqref{Eq1} in $W^{s_2,q}_0 (\Omega) \backslash \{0\}$. 
The proof of second item of Theorem \ref{MainResult}  follows from Section \ref{S3}, where, bifurcation occurs in $0$ from $\lambda_1(s_2,q)$ when we assume that $0<s_2<s_1<1<q<p<\infty$ and at infinity from $\lambda_1(s_2,q)$ when we assume that $0<s_1<s_2<1<p<q<\infty$.
The proof of the last item of Theorem \ref{MainResult} follows by combining the proof of Theorem \ref{sequence1} and Theorem \ref{sequence2}. \\

The paper is organised as follows. In Section \ref{Preli}, we introduce  the  fractional Sobolev spaces in which we will work, review some properties and embedding results, and  recall some properties related to  nonlinear eigenvalue problem invoving the $q$-Laplacian. In Section \ref{Spec}, we provide a description of the set of eigenvalues of problem \eqref{Eq1}, imposing the $L^q$-normalization  $\|u_k\|^q_q = \rho$. In Section \ref{S3}, we discuss bifurcation phenomena for problem \eqref{Eq1} and  Section \ref{Multi} is dedicated to the proof   of  multiplicity result for problem \eqref{Eq1}.

\section{Preliminaries}\label{Preli}

In this section, we introduce the definitions and some underlying properties
of the fractional Sobolev spaces, We refer the reader to \cite{di12} for further references and details on fractional Sobolev spaces. We also discuss the nonlinear eigenvalue problem for the fractional $q$-Laplace operator in an open bounded set with Lipschitz boundary.

Let $\Omega$ be an open bounded set of $\R^N, ~N\ge 1.$ Here and in the following, we identify the space $L^p(\Omega)$, $p\in [1,\infty]$ with the space of functions $u\in L^p(\R^N)$ with $u\equiv 0$ on $\R^N\setminus\Omega$.
We will work in the fractional Sobolev space $W^{s,r}(\Omega)$.  We  recall that for $0<s<1$ and $r\in[1,+\infty)$, the fractional Sobolev space $W^{s,r}(\Omega)$ is defined as
$$W^{s,r}(\Omega)=\big\{u\in L^r(\Omega):~~ [u]_{s,r,\Omega}<\infty\big\},~~$$ 
\text{where}
$$[u]_{s,r,\Omega}:=\left(\int_{\Omega}\int_{\Omega}\frac{|u(x)-u(y)|^r}{|x-y|^{N+sr}}\ dxdy\right)^{\frac{1}{r}},$$ 
is the so  called  fractional  Gagliardo seminorm. We adopte the notation $[u]_{s,r}:= [u]_{s,r,\R^N}$ when $\Omega=\R^N$. The fractional Sobolev space $W^{s,r}(\Omega)$   endowed with the norm 
$$\|u\|_{s,r}:=\left(\|u\|_r^r+[u]^r_{s,r}\right)^{1/r},$$ 
  is a reflexive Banach space. We also  define the following closed subspace of $W^{s,r}(\Omega),$
$$W^{s,r}_0(\Omega):=\big\{u\in W^{s,r}(\R^N):~~ u=0~~\text{in}~~\R^N\backslash\Omega\big\},$$ 
which is a separable Banach space endowed with the seminorm $[\cdot]_{s,r},$ which is actually a norm in $W^{s,r}_0(\Omega)$, due to  the fractional Poincar\'e inequality, 
\[
\|u\|^r_{L^r(\Omega)}\le C\int_{\R^N}\int_{\R^N}\frac{|u(x)-u(y)|^r}{|x-y|^{N+sr}}\ dxdy,
\]
with $C:=C(N,s,\Omega)$, depending on $N,s$ and $\Omega$. We  denote by $W^{-s,r'}(\Omega)$   the dual space of $W^{s,r}_0(\Omega)$, with $r$ and $r'$ satisfying the relation $\frac{1}{r}+\frac{1}{r'}=1$. We will use the following notation throughout the paper. For $u,v\in W^{s,r}(\Omega)$, we  introduce the functional $(u,v)\mapsto \cE_{r,s}(u,v)$, which  is the bilinear form associated to $(-\Delta)^s_r$~  defined by
$$\cE_{r,s}(u,v):=\langle (-\Delta)^s_ru,v\rangle_{L^2(\Omega)}= \int_{\R^N}\int_{\R^N}\frac{|u(x)-u(y)|^{r-2}(u(x)-u(y))(v(x)-v(y))}{|x-y|^{N+sr}}\ dxdy.$$
The following lemma will be used many times in this paper, see \cite[Lemma 2.1 and Lemma 2.2]{bhakta19}.
\begin{lemma}\label{Embedding} 
  Let $\Omega\subset \R^N$ be a bounded set,  $1\le q\le p\le \infty$ and $0<s_2<s_1<1$. Then, it holds that 
  \[
  W^{s_1,p}(\Omega)\subset  W^{s_2,q}(\Omega)\qquad \text{ and }\qquad W_0^{s_1,p}(\Omega)\subset  W_0^{s_2,q}(\Omega).
  \]
Moreover, there exists a constant $C:=C(N, s_1, s_2, p,q, \Omega)>0$ such that 
\[
\|u\|_{s_2,q}\le C\|u\|_{s_1,p}\quad \text{for all } \quad u\in W^{s_1,p}_0(\Omega).
%\todo{ce n'est pas un espace qu'il faut mettre}.
\] 
\end{lemma}
\noindent It follows from Lemma \ref{Embedding} that, for  $0<s_2<s_1<1<q<p<\infty$  or   $0<s_1<s_2<1<p<q<\infty$, the  embedding 
\[\text{
$W^{s_1,p}_0(\Omega) \hookrightarrow W^{s_2,q}_0(\Omega)$\qquad or \qquad $W^{s_2,q}_0(\Omega) \hookrightarrow W^{s_1,p}_0(\Omega)$}
\]
is continuous, see also  \cite[Theorem 2.1]{Chen22} and \cite[Lemma 2.1]{goel20}. 

This allows us to introduce the notion of weak solution for problem (\ref{Eq1}) in the following sense:
\begin{definition}
A function $u\in W^{s_2,q}_0(\Omega) $ or $W^{s_1,p}_0(\Omega)$ is called a weak solution of (\ref{Eq1}) if the following identity holds %for all $v\in W^{s_1,p}_0(\Omega)\cap W^{s_2,q}_0(\Omega):$
% \begin{align}\label{Eq3}\notag
%     &\iint_{\R^N\times\R^N}\frac{|u(x)-u(y)|^{p-2}(u(x)-u(y))(v(x)-v(y))}{|x-y|^{N+s_1p}}dxdy\\ 
%     &\qquad\qquad\qquad + \iint_{\R^N\times\R^N}\frac{|u(x)-u(y)|^{q-2}(u(x)-u(y))(v(x)-v(y))}{|x-y|^{N+s_2q}}dxdy=\lambda\int_{\Omega}vu~dx.
% \end{align}
\begin{equation}\label{Eq3}
    \cE_{p,s_1}(u,v)+\cE_{q,s_2}(u,v) = \lambda\int_{\Omega} |u|^{q-2}uv\ dx \qquad
\end{equation}
 \text{ for all\  $v\in W^{s_1,p}_0(\Omega)\cap W^{s_2,q}_0(\Omega)$}. 
 
 We say that $\lambda$ is an eigenvalue of problem \eqref{Eq1} in $\Omega$ if there exists a non-zero  function $u:=u_{\lambda}\in W^{s_2,q}_0(\Omega)$ or $W^{s_1,p}_0(\Omega)$ that satisfies \eqref{Eq3}. We call $u_{\lambda}$ the  eigenfunctions associated to $\lambda$ and $(\lambda,u_{\lambda})$ an eigenpair.\vspace{5pt}
\end{definition}

 We recall the following useful  inequalities.

\begin{lemma}[\cite{PL}]\label{inevec}
There exist constants $c_1,c_2$ such that for all $x_1,x_2\in\mathbb{R},$ we have the following vector inequalities for $1<r<2$
 \begin{equation*}
 (|x_2|^{r-2}x_2-|x_1|^{r-2}x_1)\cdot(x_2-x_1)\geq c_1(|x_2|+|x_1|)^{r-2}|x_2-x_1|^2,
 \end{equation*}
 and for $r>2$
 \begin{equation*}
 (|x_2|^{r-2}x_2-|x_1|^{r-2}x_1)\cdot(x_2-x_1)\geq c_2|x_2-x_1|^r.
 \end{equation*}
\end{lemma}

\begin{definition}
    An operator $\Pi: W^{s_1,p}_0(\Omega)\to W^{-s_2,p'}(\Omega)$ is said to be demi-continuous if $\Pi$ satisfies:  
    \[
\text{    $\Pi u_n\rightharpoonup\Pi u$\quad whenever\quad $u_n\in W^{s_1,p}_0(\Omega)$~~ converges to~~ $u\in W^{s_1,p}_0(\Omega),$\quad as \quad$n\to +\infty.$}
    \]
\end{definition}

Following the proof of Lemma 2.1 in \cite{Fras}, the following lemma holds.
\begin{lemma}\label{positivity}
Let $0<s_2<s_1<1<q<p<\infty$. The nonlinear operator 
\[
\text{$\Pi: W^{s_1,p}_0(\Omega)\to W^{-s_2,q'}(\Omega)\subset W^{-s_1,p'}(\Omega)$},
\]
defined by
\begin{align*}
\langle\Pi u,v\rangle=\mathcal{E}_{p,s_1}(u,v)+\mathcal{E}_{q,s_2}(u,v)
\end{align*}
 is continuous and so it is demi-continuous. In addition, the operator $\Pi$ satisfies the following condition: for any $u_n\in W^{s_1,p}_0(\Omega)$ satisfying $u_n\rightharpoonup u$ in $W^{s_1,p}_0(\Omega)$~  and~  $\limsup\limits_{n\to +\infty}\langle\Pi u_n, u_n-u\rangle\leq 0,$ 
 \[
 \text{ $u_n\to u$\quad  in \quad  $W^{s_1,p}_0(\Omega)$\quad as \quad$ n\to +\infty$.}
 \]
  The same result holds in the case $0<s_2<s_1<1<q<p<\infty$, by interchanging the role of the couples $(s_1,p)$ and $(s_2,q)$.
\end{lemma}

Finally, since we aim to show there are bifurcations branches for equation (\ref{Eq1}) emanating from the nonlinear eigenvalues of the fractional $q$-Laplacian $(-\Delta)^s_q$, we recall some spectral properties of the following eigenvalue problem 
\begin{equation}\label{dirichletlap}
	\begin{split}
	\quad\left\{\begin{aligned}
(-\Delta)^{s_2}_q u&= \lambda  |u|^{q-2}u && \text{ in \quad $\Omega$}\\
		u &=0     && \text{ on } \quad \R^N\setminus\Omega.
	\end{aligned}\right.
	\end{split}
	\end{equation}
% \begin{equation*}
% \begin{cases}
% (-\Delta)^{s_2}_q u=\lambda |u|^{q-2}u ~~~\text{in}~\Omega,\\
% u=0~~~~~~~~~~~~\text{in}~\R^N\backslash\Omega.
% \end{cases}
% \end{equation*}
It has been shown, see \cite[Proposition 2.2]{IS14} and  \cite[Theorem 1.2]{Bah} that for every $s_2\in (0,1)$ and $q\in (1,\infty)$, by means of the cohomological index and min-max theory respectively, there exists a non-decreasing sequence of variational eigenvalues of $(-\Delta)_q^{s_2}$ in $\Omega$ satisfying
\[
0<\lambda_1(s_2,q)< \lambda_2(s_2,q)\le \cdots \le \lambda_k(s_2,q)\le \cdots
\]
with $\lambda_k(s_2,q)\to \infty$ as $k\to \infty$.
% It has been shown in \cite{IS14, Bah} that there exists a nondecreasing sequence of  variational positive eigenvalues $\{\lambda_k(s_2,q)\}_k$  with  $\lambda_k(s_2,q) \to +\infty$ as $k\rightarrow +\infty$ for the following nonlinear and $q$-homogenous eigenvalue problem
% \begin{equation}\label{dirichletlap}
% 	\begin{split}
% 	\quad\left\{\begin{aligned}
% (-\Delta)^{s_2}_q u&= \lambda  |u|^{q-2}u && \text{ in \quad $\Omega$}\\
% 		u &=0     && \text{ on } \quad \R^N\setminus\Omega.
% 	\end{aligned}\right.
% 	\end{split}
% 	\end{equation}
% It is well-known that 
Moreover, the first eigenvalue $\lambda_1(s_2,q):=\lambda_{1}(s_2,q,\Omega)$ is characterized from the variational point of view by the Rayleigh quotient
\begin{equation}\label{First-q-L}
    \begin{array}{ll}
\lambda_1(s_2,q)&=\displaystyle\inf_{u\in W^{s_2,q}_0(\Omega)\backslash\{0\}}\frac{\displaystyle\int_{\R^N}\int_{\R^N}\frac{|u(x)-u(y)|^q}{|x-y|^{N+qs_2}}\ dxdy}{\displaystyle\int_{\Omega}|u(x)|^q\ dx}\vspace{10pt}\\
&=\displaystyle\inf_{\big\{u\in W^{s_2,q}_0(\Omega):~~\int_{\Omega} |u|^q\ dx=1\big\}}\int_{\R^N}\int_{\R^N}\frac{|u(x)-u(y)|^q}{|x-y|^{N+qs_2}}\ dxdy,
\end{array}
\end{equation}
and satisfies: Every eigenfunction $u_1$ corresponding to $\lambda_1(s_2,q)$ has a constant sign in $\Omega$, $\lambda_1(s_2,q)$ is simple, that is, any two eigenfunctions $u_1$ and $v_1$ corresponding to $\lambda_1(s_2,q)$ are constant multiple of each other, $u_1=\alpha v_1$, $\alpha\in\R$, see \cite[Proposition 2.1]{biswas22} . 
Furthermore,  the variational $k$-eigenvalue $\lambda_{k}(s_2,q)$ is defined by the following min-max formula, 
\begin{equation}\label{min-max}
\lambda_k(s_2,q)=\displaystyle{\inf_{A\in\Sigma_k}\sup_{u\in A}\int_{\R^N}\int_{\R^N}\frac{|u(x)-u(y)|^{q}}{|x-y|^{N+s_2q}}\ dxdy}
.\end{equation}
where $\Sigma_k$ is defined by
$$
\Sigma_k=\{A\subset\Sigma\cap D_1(s_2,q,1):~\gamma(A)\geq k\},
$$
 the set $\Sigma,$ being the class of closed symmetric  subsets of $ W^{s_2,q}_0(\Omega)\backslash\{0\},$ i.e.,
$$
\Sigma=\{A\subset W^{s_2,q}_0(\Omega)\backslash\{0\}:~~A~~\text{closed},~~ A=-A\}
$$ 
and
\[
D_{\rho}( s_2, q,\rho):=\big\{u\in W^{s_2,q}_0(\Omega)\, : \, \int_{\Omega}|u|^q\ dx=\rho\big\}.
\]
For $A\in\Sigma,$ we define    
$$\gamma(A)=\inf\{k\in\mathbb{N}~:~\exists\varphi\in C(A, \mathbb{R}^k\backslash\{0\}),~\varphi(-x)=-\varphi(x)\}.$$ If such $\gamma(A)$ does not exist, we  define $\gamma(A)=+\infty.$ The number $\gamma(A)\in\mathbb{N}\cup\{+\infty\}$ is called the {\it Krasnoselskii's genus} of $A$, see \cite{AM,ambrosio22}.
%Let us consider the family of sets 

% It is well-known  (see \cite[Proposition 2.2]{IS14}) and (see \cite[Theorem 1.2]{Bah}) that for every $s_2\in (0,1),$ by means of the cohomological index and min-max theory respectively, there exists a non-decreasing sequence of variational eigenvalue of $(-\Delta)_q^{s_2}$ in $\Omega$ satisfying
% \[
% 0<\lambda_1(s_2,q)< \lambda_2(s_2,q)\le \cdots \le \lambda_k(s_2,q)\le \cdots
% \]
%with $\lambda_k(s_2,q)\to \infty$ as $k\to \infty$. Moreover, $\lambda_1(s_2,q)$ is simple and the corresponding eigenfunction is unique up to sign and can be chosen as a positive function in $\Omega$ (\cite{LL14, F18}). 

\begin{comment}
Following the proof in \cite[Proposition 2.2]{IS14}, one shows that one has the following variational characterization of $\lambda_k(s_2,q)$,
for $k\in\mathbb{N}$,
$$
\lambda_k(s_2,q)=\displaystyle{\inf_{A\in\Sigma_k}\sup_{u\in A}\iint_{\R^N\times\R^N}\frac{|u(x)-u(y)|^{q}}{|x-y|^{N+s_2q}}dxdy}
.$$
\end{comment}

\section{The spectrum of problem (\ref{Eq1})}\label{Spec}

In this section, we provide a description of the set of eigenvalues of problem \eqref{Eq1}, imposing the $L^q$-normalization  $\|u_k\|^q_q = \rho$.
We show that equation \eqref{Eq1} has  a sequence of eigenvalues $\lambda_k(s_1, s_2, p,q;\rho)$ with $\lambda_k(s_1, s_2, p,q;\rho)\to \infty$, and with associated eigenfunctions $u_k$.
%\par \bigskip
\begin{definition}
We say that $\lambda\in\mathbb{R}$ is an eigenvalue of problem (\ref{Eq1}) if there exists a function $u_{\lambda}\in (W^{s_1,p}_0(\Omega)\cap  W^{s_2,q}_0(\Omega))\backslash\{0\}$  such that relation (\ref{Eq3}) holds.  We shall call the non-zero functions $u_{\lambda}$ an eigenfunction of problem \eqref{Eq1} associated to $\lambda.$

%\todo{Quelque chose cloche ici. On doit discuter de vive voix}\\

We say that $\lambda:=\lambda(s_1, s_2, p,q,\rho)$ is a first eigenvalue of problem \eqref{Eq1} if the corresponding eigenfunction $u_{\lambda}$ is a minimizer of the following expression for some $\rho > 0$,
\end{definition}
\begin{equation}\label{Eq4}
 \lambda_1(s_1, s_2, p,q,\rho):=\inf_{\{u\in  W^{s_1,p}_0(\Omega)\cap  W^{s_2,q}_0(\Omega), \int_\Omega |u(x)|^q\ dx = \rho\}} \frac{q}{\rho}\Big(\frac{[u]^p_{s_1,p}}{p}+\frac{[u]^q_{s_2,q}}{q} \Big). 
\end{equation}
% It holds  that $\lambda_1(s_1, s_2, p,q,1)=\lambda_1(s_2,q)$ where $\lambda_{1}(s_2,q)$ is defined in \eqref{First-q-L}. Indeed, it is clear by the definition of $\lambda_1(s_2,q)$ in \eqref{First-q-L} that $\lambda_1(s_2,q)\le \lambda_1(s_1, s_2, p,q,1)$.  Next,  for all  $t>0$ we have that 
% \begin{equation}\label{Assum}
% \lambda_1(s_1, s_2, p,q,1)\le \frac{t^{p-q}[u]^p_{s_1,p}}{p}+\frac{[u]^q_{s_2,q}}{\textcolor{red}{q}}
% \end{equation}
% Letting $t\to 0$ if $p>q$ or $t\to \infty$ if $p<q$ and applying the \textcolor{red}{infirmum over $\cC^{\infty}_c(\Omega)$~(quelle norme tu mets sur $\cC^{\infty}_c(\Omega)$??, les fonctions propres $u$ sont-elles indefiniment derivable à support compact??} in the right  hand-side of \eqref{Assum} yields $\lambda_1(s_1, s_2, p,q,1)\le \lambda_1(s_2,q)$. We conclude that $\lambda_1(s_1, s_2, p,q,1)=\lambda_1(s_2,q)$ as wanted.\\
\\
%\textbf{Voici ce que je propose comme démo de l'affirmation suivante:}\\
\\
It holds  that $\lambda_1(s_1, s_2, p,q,1)=\lambda_1(s_2,q)$ where $\lambda_{1}(s_2,q)$ is defined in \eqref{First-q-L}. Indeed we clearly have that $\lambda_1(s_1, s_2, p,q,1)\geq \lambda_1(s_2,q)$ since a positive term is added. On the other hand, consider $u=\frac1t e_1$ (where $e_1$ is the first eigenfunction of $(-\Delta)^{s_2}_q$ associated to $\lambda_1(s_2,q)$), we get
$$\lambda_1(s_1, s_2, p,q,1)\leq \frac{\frac{[e_1]^p_{s_1,p}}{pt^p}+\frac{[e_1]^q_{s_2,q}}{qt^q}}{\frac{1}{qt^q}\int_{\Omega}|e_1|^qdx}=\frac{\frac{[e_1]^p_{s_1,p}}{pt^{p-q}}+\frac{[e_1]^q_{s_2,q}}{q}}{\frac1q\int_{\Omega}|e_1|^qdx}\to \lambda_1(s_2,q) $$
as $t\to \infty$ if $p>q$ and as $t\to 0$ if $q<q.$ We conclude that $\lambda_1(s_1, s_2, p,q,1)=\lambda_1(s_2,q)$ as wanted.\\
% \\
% We note that there exists $u_{\lambda_1}\in W^{s_1,p}_0(\Omega)\cap  W^{s_2,q}_0(\Omega)$ such that  $\lambda_1(s_1, s_2, p,q;\rho)$ satisfies 
% \begin{equation*}
%     \frac{[u_{\lambda_1}]^p_{s_1,p}}{p}+\frac{[u_{\lambda_1}]^q_{s_2,q}}{q}=\lambda_1(s_1, s_2, p,q,\rho)\int_{\Omega}|u_{\lambda_1}(x)|^q\ dx.
% \end{equation*}

We have the following nonexistence result for problem \eqref{Eq1}.
\begin{proposition}[Nonexistence] Let $\lambda_1(s_2,q)$ be the first eigenvalue of 
\begin{equation*}
(-\Delta)^{s_2}_q u~= ~\lambda  |u|^{q-2}u ~~~~ \text{ in \quad $\Omega$},\quad\qquad
		u ~=0     ~~ \text{ on } \quad \R^N\setminus\Omega.
	\end{equation*}
 
If it holds that $\lambda\leq \lambda_1(s_2,q)$, then problem (\ref{Eq1}) has no nontrivial solutions.
\end{proposition}

\begin{proof}
Suppose by contradiction that there exists $\lambda<\lambda_1(s_2,q)$ which is an eigenvalue of problem (\ref{Eq1}) with $u_{\lambda}\in (W^{s_1,p}_0(\Omega)\cap  W^{s_2,q}_0(\Omega))\backslash\{0\}$ the corresponding eigenfunction. Letting $v=u_{\lambda}$ in relation (\ref{Eq3}), we then have
\begin{equation*}
    \iint_{\R^N\times\R^N}\frac{|u_{\lambda}(x)-u_{\lambda}(y)|^p}{|x-y|^{N+s_1p}}\ dxdy+\iint_{\R^N\times\R^N}\frac{|u_{\lambda}(x)-u_{\lambda}(y)|^q}{|x-y|^{N+qs_2}}\ dxdy=\lambda\int_{\Omega}|u_{\lambda}(x)|^q\ dx.
\end{equation*}
On the other hand, the Poincar\'e inequality yields
\begin{equation}\label{Eq5}
\lambda_1(s_2,q)\int_\Omega|u_{\lambda}(x)|^q\ dx\leq \iint_{\R^N\times\R^N}\frac{|u_{\lambda}(x)-u_{\lambda}(y)|^q}{|x-y|^{N+qs_2}}\ dxdy.
\end{equation}
Subtracting  both sides of (\ref{Eq5}) by $\lambda\int_{\Omega}|u_{\lambda}(x)|^q\ dx$,~ it follows that 
$$\left(\lambda_1(s_2,q)-\lambda\right)\int_\Omega|u_{\lambda}(x)|^qdx\leq \iint_{\R^N\times\R^N}\frac{|u_{\lambda}(x)-u_{\lambda}(y)|^q}{|x-y|^{N+qs_2}}\ dxdy-\lambda \int_\Omega|u_{\lambda}(x)|^q\ dx.$$ 
This implies that
\[
\begin{array}{ll}
    0<\left(\lambda_1(s_2,q)-\lambda\right)\displaystyle\int_\Omega|u_{\lambda}(x)|^q\ dx&\leq \displaystyle\iint_{\R^N\times\R^N}\frac{|u_{\lambda}(x)-u_{\lambda}(y)|^q}{|x-y|^{N+qs_2}}\ dxdy-\lambda \displaystyle\int_\Omega|u_{\lambda}(x)|^q\ dx\vspace{10pt}\\
    &\qquad+\displaystyle\iint_{\R^N\times\R^N}\frac{|u_{\lambda}(x)-u_{\lambda}(y)|^p}{|x-y|^{N+s_1p}}\ dxdy=0.
\end{array}
\]
Hence, $\lambda<\lambda_1(s_2,q)$ is not an eigenvalue of problem (\ref{Eq1}) with $u_{\lambda}\neq 0.$ The proof for the case $\lambda=\lambda_1(s_2,q)$ follows from the definition of $\lambda_1(s_2,q)$. Indeed, taking $v=u_{\lambda_{1}}$ in (\ref{Eq3}), we have
\begin{equation*}
    \iint_{\R^N\times\R^N}\frac{|u_{\lambda_1}(x)-u_{\lambda_1}(y)|^p}{|x-y|^{N+s_1p}}\ dxdy+\iint_{\R^N\times\R^N}\frac{|u_{\lambda_1}(x)-u_{\lambda_1}(y)|^q}{|x-y|^{N+qs_2}}\ dxdy=\lambda_1(s_2,q)\int_{\Omega}|u_{\lambda_1}(x)|^q\ dx.
\end{equation*}
It follows from the definition of $\lambda_1(s_2,q)$ that
\begin{equation*}
    \iint_{\R^N\times\R^N}\frac{|u_{\lambda_1}(x)-u_{\lambda_1}(y)|^p}{|x-y|^{N+s_1p}}\ dxdy\le \lambda_1\int_{\Omega}|u_{\lambda_1}(x)|^q\ dx-\lambda_1\int_{\Omega}|u_{\lambda_1}(x)|^q\ dx=0.
\end{equation*}
It follows from Poincar\'e inequality that
\begin{equation}
0\le \int_\Omega|u_{\lambda_1}(x)|^p\ dx\leq C\iint_{\R^N\times\R^N}\frac{|u_{\lambda_1}(x)-u_{\lambda_1}(y)|^p}{|x-y|^{N+ps_1}}\ dxdy =0\implies u_{\lambda_1}\equiv 0\quad \text{ in }\ \Omega.
\end{equation}
This shows that any eigenvalue of problem (\ref{Eq1}) satisfies $\lambda \in (\lambda_1(s_2,q),\infty)$. 
\end{proof}

The Palais-Smale condition plays an important role in the minimax argument, and we recall here its definition.
\begin{definition}
A $C^1$ functional $I$ defined on a smooth submanifold $M$ of a Banach space $X$ is said to satisfy the Palais-Smale condition on $M$ if any sequence $\{u_n\}\subset M$ satisfying that $\{I(u_n)\}_n$ is bounded and $\big(I\big|_M\big)'(u_n)\rightarrow 0$ as $n\rightarrow +\infty$ has a convergent subsequence.
\end{definition}

Next, we start the discussion about the existence of eigenvalues for problem (\ref{Eq1}). We note that these eigenvalues depend on $\rho$, from the $L^q$-normalization  $\int_{\Omega}|u(x)|^q\ dx=\rho.$ We define the  energy functional $F_{\lambda}: W^{s_1,p}_0(\Omega)\cap  W^{s_2,q}_0(\Omega)\rightarrow \mathbb{R}$ associated to relation (\ref{e2}) by
\begin{equation}\label{functional}
    F_{\lambda}(u)=\frac{1}{p}[u]^p_{s_1,p}+\frac{1}{q}[u]^q_{s_2,q}-\frac{\lambda}{q}\int_{\Omega}|u|^qdx.
\end{equation}
A standard arguments can be used to show that $F_{\lambda}\in C^1(W^{s_1,p}_0(\Omega),\mathbb{R})$ with its derivative given by $$\langle F'_{\lambda}(u),v\rangle=\cE_{p,s_1}(u,v)+\cE_{q,s_2}(u,v)-\lambda\int_{\Omega}|u|^{q-2}u~v~dx,$$ for all $v\in W^{s_1,p}_0(\Omega)\cap W^{s_2,q}_0(\Omega)$. Thus, we note that $\lambda$ is an eigenvalue of problem (\ref{Eq1}) if and only if $F_{\lambda}$ possesses a nontrivial critical point. 

We further split the discussion into two cases, whether, $0<s_2< s_1<1<q<p<\infty$ ~ or ~  $0<s_1< s_2<1<p<q<\infty$.

\subsection{The case: $0<s_2< s_1<1<q<p<\infty$ }

In this case one  show that for each $\lambda>0,$ the functional $F_{\lambda}$ defined in (\ref{functional}) is coercive. 

%In this section we prove that every $\lambda > \lambda_1(s_2,q)$ is a first eigenvalue of problem \eqref{Eq1}.

\begin{lemma}\label{coercivite}
 For each $\lambda>0,$ the functional $F_{\lambda}$ defined in (\ref{functional}) is coercive.
\end{lemma}
\begin{proof}
Using the fact that $W_0^{s_1,p}(\Omega)\hookrightarrow L^q(\Omega)$ (see Lemma \ref{Embedding}), we have
\begin{align*}
   F_{\lambda}(u)\geq \frac{[u]^p_{s_1,p}}{p}-\frac{\lambda}{q}\int_{\Omega}|u|^qdx&\geq \frac{[u]^p_{s_1,p}}{p}-\lambda C\frac{[u]^q_{s_1,q}}{q}.
\end{align*}
Therefore, $F_{\lambda}(u)\to +\infty$ as $[u]_{s_1,p}\to +\infty$ since $p>q.$
\end{proof}
\begin{remark}
Note that the functional $F_{\lambda}$ is not bounded below if $p<q$ and $\lambda>\lambda_1(s_2,q)$. Indeed,  for  $u=u_1,$ the first eigenfunction of problem (\ref{dirichletlap}) with $\int_{\Omega}|u_1|^qdx=1,$ we have 
$$F_{\lambda}(tu_1)=\frac{t^p}{p}[u_1]^p_{s_1,p}+\frac{t^q}{q}(\lambda_1(s_2,q)-\lambda)\rightarrow -\infty\qquad \text{as \quad $t\rightarrow+\infty.$}$$
This case will be treated in Section \ref{case2}, on a subset of $W^{s_1,p}_0(\Omega)\cap  W^{s_2,q}_0(\Omega)$, the so-called Nehari manifold, since we cannot apply the Direct Method in the Calculus of Variations in order to find critical points for the functional $F_{\lambda}$. 
\end{remark}

\begin{theorem}\label{Existence1}
Every $\lambda\in  (\lambda_1(s_2,q),+\infty)$ is an eigenvalue of problem (\ref{Eq1}).
\end{theorem}
\begin{proof}
As already  mentioned above, standard arguments show that $F_{\lambda}\in C^1(W^{s_1,p}_0(\Omega),\mathbb{R})$ with its derivative given by 
$$\langle F'_{\lambda}(u),v\rangle=\cE_{p,s_1}(u,v)+\cE_{q,s_2}(u,v)-\lambda\int_{\Omega}|u|^{q-2}u~v~dx,$$
for all $v\in W^{s_1,p}_0(\Omega)\subset W^{s_2,q}_0(\Omega).$ On the other hand $F_{\lambda}$ is weakly lower semi-continuous on $W^{s_1,p}_0(\Omega)\subset W^{s_2,q}_0(\Omega)$ since $F_{\lambda}$ is a continuous convex functional. This fact and Lemma \ref{coercivite} allow one to apply the Direct Method of Calculus of Variations   to obtain the existence of global minimum point of $F_{\lambda}$. We denote by $u_0$ such a global minimum point, i.e, 
$$F_{\lambda}(u_0)=\min\limits_{u\in W^{s_1,p}_0(\Omega) }F_{\lambda}(u).$$ 
We observe that for $u_0=tw_1$, where $w_1$ stands for the $L^q$-normalized associated eigenfunction of $\lambda_1(s_2,q)$, we have, 
$$F_{\lambda}(u_0)=F_{\lambda}(tw_1)=\frac{t^p}{p}[w_1]^p_{s_1,p}+\frac{t^q}{q}(\lambda_1(s_2,q)-\lambda)<0$$ for $t$ small enough. So there exists $u_{\lambda}\in W^{s_1,p}_0(\Omega)$ such that $F_{\lambda}(u_{\lambda})<0.$ But $F_{\lambda}(u_0)\leq F_{\lambda}(u_{\lambda})<0,$ which implies that $u_0\in W^{s_1,p}_0(\Omega)\backslash\{0\}.$ We also have that 
$$\text{$\langle F'_{\lambda}(u_0),v\rangle=0$,\quad for all\quad \ $v \in W^{s_1,p}_0(\Omega).$}$$ 
This concludes the proof of Theorem \ref{Existence1}.
\end{proof}

\begin{proposition}\label{constantsign}
     Every eigenfunction $u_{\lambda}$ associated to $\lambda\in (\lambda_1(s_2,q),\infty)$ is positive or negative in $\Omega$.
\end{proposition}
\begin{proof}
  Let $u_{\lambda}\in W^{s_1,p}_0(\Omega)\setminus\{0\}$ be an eigenfunction associated to $\lambda\in(\lambda_1(s_2,q),\infty). $ Then,
  $$ \frac{1}{p}\iint_{\R^N\times\R^N}\frac{|u_{\lambda}(x)-u_{\lambda}(y)|^p}{|x-y|^{N+s_1p}}\ dxdy+\frac{1}{q}\iint_{\R^N\times\R^N}\frac{|u_{\lambda}(x)-u_{\lambda}(y)|^q}{|x-y|^{N+s_2q}}\ dxdy=\frac{\lambda}{q}\int_{\Omega}|u_{\lambda}(x)|^{q}\ dx,$$ which means that $u_{\lambda}$ achieves the infimum for $F_{\lambda}$ with $\lambda:=\lambda(s_1, s_2, p,q)$  in Theorem \ref{Existence1}. From the triangle inequality, we have that 
  $$||u(x)|-|u(y)||\le |u(x)-u(y)|\quad \text{ for all }\quad u\in W^{s_1,p}_0(\Omega)\subset W^{s_2,q}_0(\Omega),$$ 
  which implies that
     $$[|u_{\lambda}|]^p_{s_1,p}\le[u_{\lambda}]^p_{s_1,p}\quad\text{ and }\quad [|u_{\lambda}|]^q_{s_2,q}\le [u_{\lambda}]^q_{s_2,q}$$
   and thus  $|u_{\lambda}|\in W^{s_1,p}_0(\Omega)\subset W^{s_2,q}_0(\Omega)$. Hence, it follows that $F_{\lambda}(|u_{\lambda}|)\le F_{\lambda}(u_{\lambda})$. This yields,
 $$  
   F_{\lambda}(|u_{\lambda}|)\le F_{\lambda}(u_{\lambda})=\min\limits_{v\in W^{s_1,p}_0(\Omega) }F_{\lambda}(v)\le F_{\lambda}(|u_{\lambda}|),
  $$ 
%   By the definition of $\lambda$, we have that 
   % \begin{equation}
   %  \lambda_1(s_1, s_2, p,q;\rho)= \frac{[|u_{\lambda}|]^p_{s_1,p}}{p}+\frac{[|u_{\lambda}|]^q_{s_2,q}}{q},~~\textbf{il manque le inf dans l'expression}
   % \end{equation}
   showing that $|u_{\lambda}|$ also achieves the infimum for  $F_{\lambda}$ and hence, $u_{\lambda}$ must has a  constant sign in $\Omega$. 
\end{proof}

The proofs of the following two theorems rely on \cite[Proposition 10.8]{AM}.

\begin{theorem}\label{sequence1}
There exists a nondecreasing sequence of critical values $c_k(s_1,s_2,p,q;\rho)$ with associated nonlinear eigenvalues $\lambda_k(s_1,s_2,p,q;\rho)\to +\infty,$ as $k\to +\infty$ and with corresponding eigenfunctions $u_k\in W^{s_1,p}_0(\Omega)$ with  $\int_{\Omega}|u_k|^p\ dx=\rho,$ for problem (\ref{Eq1}).
\end{theorem}

\begin{proof} We define:
\begin{align*}
    D_{\rho}(s_1, s_2, p,q)&=\big\{u\in W^{s_1,p}_0(\Omega)\, : \, \int_{\Omega}|u|^q\ dx=\rho\big\},\vspace{5pt}\\
    \Sigma_k(s_1,s_2,p,q)&=\big\{A \subset D_\rho(s_1,s_2,p,q), \quad A \in \Sigma\quad \hbox{ and} \quad \gamma(A) \geq k\big\},
\end{align*}
    where $\Sigma=\{A\subset W^{s_1,p}_0(\Omega): \ A\,\text{closed},\quad A=-A\}$ and set 
\begin{equation}\label{ck}
c_k(s_1,s_2,p,q;\rho)=\inf_{A\in\Sigma_k(s_1,s_2,p,q)}\sup_{u\in A}\Big(\frac{[u]^p_{s_1,p}}{p}+\frac{[u]^q_{s_2,q}}{q} \Big)>0.
\end{equation}
Let us show that 
$$
I(u)=\frac{[u]^p_{s_1,p}}{p}+\frac{[u]^q_{s_2,q}}{q}
$$ 
satisfies the Palais-Smale (PS) condition on $D_{\rho}(s_1,s_2,p,q)$. 
Let $\{u_n\}\subset D_{\rho}(s_1,s_2,p,q)$ be a (PS) sequence, i.e, for all $n,$ there exists a constant $K>0$ such that
$$|I(u_n)|\leq K$$
and 
$$
\text{$(I\big|_{D_\rho})'(u_n)\rightarrow 0$\quad  in\quad  $\Big(W^{s_1,p}_0(\Omega)\Big)'\equiv W^{-s_1,p'}(\Omega)$ ~~ as ~~$n \rightarrow \infty.$}
$$
First, we show that $\{u_n\} \subset D_{\rho}(s_1,s_2,p,q)$ is bounded in $W^{s_1,p}_0(\Omega).$ Since $u_n\in W^{s_1,p}_0(\Omega)\subset W^{s_2,q}_0(\Omega),$ with the Poincar\'e inequality, we have $\rho=\int_{\Omega}|u_n|^q\ dx\leq C[u_n]^q_{s_2,q} $ and it follows that 
$$
K\geq |I(u_n)|\geq \frac{q}{p}[u_n]^p_{s_1,p}+\frac{1}{C}\int_{\Omega}|u_n|^qdx=\frac{q}{p}[u_n]^p_{s_1,p}+\frac{\rho}{C}.
$$ 
Then $\{u_n\} \subset D_{\rho}(s_1,s_2,p,q)$ is bounded in $ W^{s_1,p}_0(\Omega).$ We can assume that up to a subsequence, still denoted $\{u_n\}$, there exists $u\in W^{s_1,p}_0(\Omega)$ such that $u_n \rightharpoonup u$ in $W^{s_1,p}_0(\Omega)$ and $u_n \to u$ in $L^q(\Omega)$. Now, we show that $u_n$ converges strongly to $u$ in $W^{s_1,p}_0(\Omega).$ Since 
$(I\big|_{D_\rho})'(u_n)\rightarrow 0$ in $W^{-s_1,p'}(\Omega)$ as $n\rightarrow +\infty,$ there exists $\mu_n \in \mathbb R$ and $\varepsilon_n \to 0$ in $W^{-s_1,p'}_0(\Omega)$ such that 
\[
\text{ $I'(u_n)v - \mu_n \int_\Omega |u_n|^{q-2}u_n v\ dx = \langle \varepsilon_n,v\rangle$\quad for all\quad $v\in W_0^{s_2,q}(\Omega)$.  }
\]
We have $I'(u_n)u_n - \mu_n\int_\Omega |u_n|^q \to 0$, and since $I'(u_n)u_n \le C_1 I(u_n)\le C_2$  it follows that $|\mu_n| \le C$. From this we obtain that $I'(u_n)(u_n-u)\rightarrow 0$ and $I'(u)(u_n-u)\rightarrow 0$ as $n\rightarrow+\infty.$ 
Therefore, using Lemma \ref{inevec}, we have:

If $1<q<p<2,$ then,
\[
\begin{split}
 o(1) &=   \langle I'(u_n)-I'(u), u_n-u\rangle\\
  &\ge \langle (-\Delta)_p^{s_1}u_n-(-\Delta)_p^{s_1}u,u_n-u\rangle  +\langle (-\Delta)_q^{s_2}u_n-(-\Delta)_q^{s_2}u,u_n-u\rangle   \\
  &\ge \|u_n-u\|_{s_1,p}\left( \|u_n\|_{s_1,p}+\|u\|_{s_2,p}\right)^\frac{p-2}{p}+\|u_n-u\|_{s_2,q}\left( \|u_n\|_{s_2,q}+\|u\|_{s_2,q}\right)^\frac{q-2}{p}\\
  &\ge \|u_n-u\|_{s_1,p}\left( \|u_n\|_{s_1,p}+\|u\|_{s_1,p}\right)^\frac{p-2}{p}
\end{split}
\]
so that 
\[
\begin{array}{ll}
\|u_n-u\|_{s_1,p}&\le C\left( \|u_n\|_{s_1,p}+\|u\|_{s_1,p}\right)^\frac{2-p}{p}\langle I'(u_n)-I'(u), u_n-u\rangle\vspace{10pt}\\
&\le C \langle I'(u_n)-I'(u), u_n-u\rangle \\
&=~o(1).
\end{array}
\]

If $p>q\ge 2,$ still by Lemma \ref{inevec}, we have
\[
\begin{array}{ll}
    \|u_n-u\|_{{s_1,p}}&\le \|u_n-u\|_{{s_1,p}}+\|u_n-u\|_{{s_2,q}}\vspace{10pt}\\
    &\le \langle (-\Delta)_p^{s_1}u_n-(-\Delta)_p^{s_1}u,u_n-u\rangle  +\langle (-\Delta)_q^{s_2}u_n-(-\Delta)_q^{s_2}u,u_n-u\rangle \vspace{10pt}\\
    &\le C \langle I'(u_n)-I'(u), u_n-u\rangle\vspace{10pt}\\
&= o(1)
\end{array}
\]
Therefore, it holds for every $1<q<p<\infty$ and $0<s_2< s_1<1$, that 
$$o(1)=\langle I'(u_n)-I'(u), u_n-u\rangle\geq c_2\|u_n-u\|^p_{s_1,p}.$$ This shows that $u_n$ converges strongly to $u$ in $W^{s_1,p}_0(\Omega)$ as $n\rightarrow+\infty. $ %since $\langle I'(u_n)-I'(u), u_n-u\rangle\rightarrow 0$ as $n\rightarrow+\infty.$\\

 In order to complete the proof, let us show that if $c=c_k(s_1,s_2,p,q;\rho)=\dots=c_{k+m-1}(s_1,s_2,p,q;\rho),$ then the set $K_c$ of critical points of $I$ at the critical level $c$ has a genus $\gamma(K_c)\geq m.$
We consider the level set at $c$, $$K_c:=\{u\in D_{\rho}(s_1,s_2,p,q):~I(u)=c~,~I'(u)=0\}.
$$ 
We have that $K_c$ is compact since the functional $I$ satisfies the Palais-Smale condition and $0\notin K_c$ since $c>0=I(0)$. In addition, we have $I(u)=I(-u)$. Hence $K_c\in \Sigma$. 

Assume by contradiction that $\gamma(K_c) \le  m-1$. Take $A_\varepsilon \in \Sigma_{k+m-1}$ such that $\sup_{A_\varepsilon} I(u) \le c+\varepsilon$. By the properties of the genus, there exists a $\delta$-neighborhood $N_\delta$ of $K_c$ such that $\gamma(N_\delta) = \gamma(K_c)$, and
$\gamma(A_\varepsilon \setminus N_\delta) \ge \gamma(A_\varepsilon) - \gamma(N_\delta) \ge k+m-1 - (m-1) = k$. By the deformation theorem there exists a homeomorphism $\eta(1,\cdot)$ such that $I(u) \le c-\varepsilon$, for $u \in \eta(1,A_\varepsilon \setminus N_\delta)$. Then we arrive at the contradiction
$$
c = \inf_{A\in \Sigma_k}\sup_{u \in A} I(u)\le \sup_{\eta(1,A_\varepsilon \setminus N_\delta)}I(u) \le c-\varepsilon
.$$
Hence, $\gamma(K_c) \ge m$.
\par
With a compactness argument one shows that $c_k(s_1,s_2,p,q;\rho) \to +\infty$ as $k \to +\infty$. 
%\par \smallskip
For the corresponding eigenvalues $\lambda_k(s_1,s_2,p,q;\rho)$ we then have
$$ 
[u_n]^p_{s_1,p} + [u_n]^q_{s_2,q}
= \lambda_k(s_1,s_2,p,q;\rho) \int_\Omega |u_n|^q\ dx = \lambda_k(s_1,s_2,p,q;\rho) \, \rho.
$$
Thus $\lambda_k(s_1,s_2,p,q;\rho) \, \rho > c_k(s_1,s_2,p,q;\rho)$, for all $k$ (and fixed $\rho$), and hence also $\lambda_k(s_1,s_2,p,q;\rho) \to +\infty$ as $k \to +\infty$.
\end{proof}

\subsection{The case: $0<s_1< s_2<1<p<q<\infty$ }\label{case2}

In this case, we do not have  coercivity on $W^{s_2,q}_0(\Omega)\setminus\{0\}$ for functional $J_{\lambda}$ eventhough it belongs to $\cC^1(W^{s_2,q}_0(\Omega)\setminus\{0\};\R)$. To prove that $F_{\lambda
}$ has a critical point in $W^{s_2,q}_0(\Omega)\setminus\{0\}$, we constrain $F_{\lambda}$ on the Nehari set
\begin{equation*}
\begin{split}
    \mathcal{N}_{\lambda}&=\{u\in W^{s_2,q}_0(\Omega)\setminus\{0\},~\langle F'_{\lambda}(u),u\rangle=0\}\\
    &=\{u\in W^{s_2,q}_0(\Omega)\setminus\{ 0\},~[u]^p_{s_1,p}+[u]^q_{s_2,q}=\lambda\int_{\Omega}|u|^q\ dx\}.
\end{split}
\end{equation*}
Note that on $\mathcal{N}_{\lambda},$ the functional $F_{\lambda}$ reads as $$F_{\lambda}(u)=(\frac{1}{p}-\frac{1}{q})[u]^p_{s_1,p}>0.$$
This shows at once that $F_{\lambda}$ is coercive in the sense that if $u\in\mathcal{N}_{\lambda}$ satisfies $[u]_{s_1,p}\rightarrow+\infty,$ then $F_{\lambda}(u)\rightarrow +\infty.$

Next, we define the quantity
$$\mathfrak{M}=\inf\limits_{u\in\mathcal{N}_{\lambda}}F_{\lambda}(u)$$
and  show through a series of propositions that~ $\mathfrak{M} $ is attained by a function $u\in\mathcal{N}_{\lambda}$, which is a critical point of $F_{\lambda}$ considered on the whole space $W^{s_2,q}_0(\Omega)\subset W^{s_1,p}_0(\Omega)$ and in fact is a solution of (\ref{Eq1}).

\begin{proposition}\label{nonvide}
The set $\mathcal{N}_{\lambda}$ is not empty for any $\lambda>\lambda_1(s_2,q).$
\end{proposition}
\begin{proof}
Since $\lambda>\lambda_1(s_2,q)$ there exists $u\in W^{s_2,q}_0(\Omega)$ not identically zero such that $[u]^q_{s_2,q}<\lambda\int_{\Omega}|u|^qdx.$ We then see that $tu\in \mathcal{N}_{\lambda}$ for some $t>0.$ Indeed, $tu\in \mathcal{N}_{\lambda}$ is equivalent to
$$t^p[u]^p_{s_1,p}+t^q[u]^q_{s_2,q}=t^q\lambda\int_{\Omega}|u|^q\ dx,$$ which is solved for $t=\left(\frac{[u]^p_{s_1,p}}{\lambda\int_{\Omega}|u|^q\ dx-[u]^q_{s_2,q}}\right)^{\frac{1}{q-p}}>0.$
\end{proof}

\begin{proposition}\label{minborne}
Every minimizing sequence for $F_{\lambda}$ on $\mathcal{N}_{\lambda}$ is bounded in $W^{s_2,q}_0(\Omega).$
\end{proposition}
\begin{proof}
Let $\{u_n\}_{n\geq0}\subset \mathcal{N}_{\lambda}$ be a minimizing sequence of $F_{\lambda}|_{\mathcal{N}_{\lambda}}$, i.e. $F_{\lambda}(u_n)\rightarrow \mathfrak{M}=\displaystyle{\inf_{v\in\mathcal{N}_{\lambda}}}F_{\lambda}(v).$
Then
\begin{equation}\label{m1}
\lambda\int_{\Omega}|u_n|^q~dx-
[u_n]^q_{s_2,q}=[u_n]^p_{s_1,p}\rightarrow\left(\frac{1}{p}-\frac{1}{q}\right)^{-1}\mathfrak{M},~\text{as $n\rightarrow +\infty$}.
\end{equation}
Suppose on the contrary that $\{u_n\}_{n\geq0}$ is not bounded in $W^{s_2,q}_0(\Omega)$ i.e. $[u_n]^q_{s_2,q}\rightarrow +\infty$ as $n\rightarrow +\infty$. Then we have $\displaystyle{\int_{\Omega}}|u_n|^q\ dx
\rightarrow\infty$ as $n\rightarrow +\infty$, using relation (\ref{m1}). Next, we set $w_n=\frac{u_n}{\|u_n\|_q}.$ Since $[u_n]^q_{s_2,q}<\lambda\displaystyle{\int_{\Omega}}|u_n|^q~dx
$, we deduce that $[w_n]^q_{s_2,q}<\lambda,$ for each $n$. Hence $\{w_n\}\subset W^{s_2,q}_0(\Omega)$ is uniformly bounded in $W^{s_2,q}_0(\Omega).$ Therefore there exists $w_0\in W^{s_2,q}_0(\Omega)$ such that passing to a subsequence that we still label $\{w_n\}_n$, 
\[
\text{
$w_n\rightharpoonup w_0$\quad in \quad $W^{s_2,q}_0(\Omega)\subset W^{s_1,p}_0(\Omega)$ \quad and \quad  $w_n\rightarrow w_0$ in $L^q(\Omega)$ }
\]
thanks to the compactness of the  Sobolev embedding. Since $p<q$, by  H\"older inequality, we have that
$
\|w_n\|_{L^p(\Omega)}\le |\Omega|^{\frac{q-p}{pq}}\|w_n\|_{L^q(\Omega)}\le C 
$
and $w_n\to w_0$ in $L^p(\Omega)$.
Furthermore, as we know that $\lambda\int_{\Omega}|u_n|^q\ dx-
[u_n]^q_{s_2,q}$ is bounded as $n\to \infty$ and $p<q$, it follows from (\ref{m1}) combined with the Poincar\'e inequality that
\[
\begin{split}
0\le\int_{\Omega}|w_n|^p\ dx&=  \frac{1}{\|u_n\|^p_{L^q(\Omega)}}\int_{\Omega}|u_n|^p\ dx\\
&\le \frac{C}{\|u_n\|^p_{L^q(\Omega)}}[u_n]^p_{s_1,p}=C\ \frac{\displaystyle\lambda\int_{\Omega}|u_n|^q \ dx-
[u_n]^q_{s_2,q}}{\|u_n\|^p_{L^q(\Omega)}}\to 0\quad \text{ as }\quad n\to \infty.
\end{split}
\]
Therefore,  $w_n \to 0$ in $L^p(\Omega)$ and  consequently $w_0=0$. This contradicts the normalization of $\|w_n\|_{L^q(\Omega)}=1$. Hence the sequence $\{u_n\}_n$ is bounded in $W^{s_2,q}_0(\Omega)$. This completes the proof of Proposition \ref{minborne}.
\end{proof}
%\todo{à continuer}

\begin{proposition}\label{minpositive}
We have ~ $\mathfrak{M}=\inf\limits_{u\in\mathcal{N}_{\lambda}}F_{\lambda}(u)>0.$
\end{proposition}
\begin{proof}
Assume by contradiction that $\mathfrak{M}=0$. Then, with $\{u_n\}_{n\geq 0}$ as in Proposition \ref{minborne}, we have
\begin{equation}\label{e6}
0<\lambda\int_{\Omega}|u_n|^q \ dx-
[u_n]^q_{s_2,q}=[u_n]^p_{s_1,p}~\rightarrow ~0, \quad\text{as ~~~$n\rightarrow +\infty$}.
\end{equation}
By Proposition \ref{minborne}, we know that $\{u_n\}_{n\geq 0}$ is bounded in $W^{s_2,q}_0(\Omega).$ Therefore there exists $u_0\in W^{s_2,q}_0(\Omega)$ such that up to a subsequence,  $u_n \rightharpoonup u_0$ in $W^{s_2,q}_0(\Omega)\subset W^{s_1,p}_0(\Omega)$ and $u_n\rightarrow u_0$ in $L^q(\Omega).$
Thus, by the weakly lower semicontinuity of the norm,
$$\|u_0\|^p_{L^p(\Omega)}\le C[u_0]^p_{s_1,p}\leq C\liminf\limits_{n\rightarrow +\infty}[u_n]^p_{s_1,p}=0.$$ Consequently  $u_0=0$. Passing to a subsequence again, we have that $u_n \rightharpoonup 0$ in $W^{s_2,q}_0(\Omega)\subset W^{s_1,p}_0(\Omega)$ and $u_n\rightarrow 0$ in $L^q(\Omega) $ and also $u_n\rightarrow 0$ in $L^p(\Omega) $ since $p<q.$  Next,  writing again $w_n=\frac{u_n}{\|u_n\|_q}$, we have that $\|w_n\|_{q}=1$ for all $n$ and as in Proposition \ref{minborne}, $\{w_n\}$ is uniformly bounded in $W^{s_2,q}_0(\Omega).$ Hence $w_n \rightharpoonup w_0$ in $W^{s_2,q}_0(\Omega)$ and $w_n\rightarrow w_0$ in $L^q(\Omega)$ and $w_n\rightarrow w_0$ in $L^p(\Omega)$ since $p<q$.  By Poincar\'e inequality, we have
\[
\begin{split}
    0\le  \|w_n\|_p &\le C  [w_n]^p_{s_1,p}\\
 %   & =C |\Omega|^{\frac{p-q}{q}}\left(\lambda\int_{\Omega}|w_n|^q~dx-
%[w_n]^q_{s_2,q}  \right)\\
& =\frac{C }{\|u_n\|^p_q}\left(\lambda\int_{\Omega}|u_n|^q~dx-
[u_n]^q_{s_2,q}  \right) \\
&=C \|u_n\|^{q-p}_q\left(\lambda-
[w_n]^q_{s_2,q}  \right) \to 0\quad \text{ as }\quad n\to \infty.
\end{split}
\]
We deduce  taking the limit that $w_0=0$, which contradicts the normalization $\|w_n\|_{q}=1$. This completes the proof of Proposition \ref{minpositive}.
\end{proof}

\begin{proposition}\label{minatteind}
 There exists $u\in \mathcal{N}_{\lambda}$ such that
$F_{\lambda}(u)=\mathfrak{M}.$
\end{proposition}
\begin{proof}
Let $\{u_n\}_{n\geq 0}\subset\mathcal{N}_{\lambda}$ be a minimizing sequence, i.e., $F_{\lambda}(u_n)\rightarrow \mathfrak{M}$ as $n\rightarrow\infty.$ Thanks to Proposition \ref{minborne}, we have that $\{u_n\}$ is bounded in $W_0^{s_2,q}(\Omega).$ It follows that there exists $u_0\in W_0^{s_2,q}(\Omega)$ such that $u_n\rightharpoonup u_0$ in $W_0^{s_2,q}(\Omega)\subset W_0^{s_1,p}(\Omega)$ and strongly in $L^q(\Omega).$ The results in the two propositions above guarantee that $F_{\lambda}(u_0)\leq \displaystyle{\lim_{n\rightarrow\infty}}\inf F_{\lambda}(u_n)=\mathfrak{M}.$ Since for each $n$ we have $u_n\in\mathcal{N}_{\lambda}$, then
\begin{equation}\label{att1}
[u_n]^p_{s_1,p}+[u_n]^q_{s_2,q}=\lambda \int_{\Omega}|u_n|^q~dx~~~\text{for all $n$.}
\end{equation}
Assuming $u_0\equiv 0$ on $\Omega$ implies that $ \displaystyle{\int_{\Omega}}|u_n|^q~dx\rightarrow 0$ as $n\rightarrow +\infty$, and by relation $(\ref{att1})$ we obtain that $[u_n]^q_{s_2,q}\rightarrow 0$ as $n\rightarrow +\infty.$ Combining this with the fact that $u_n \rightharpoonup 0$ in $W_0^{s_2,q}(\Omega)$, we deduce that $u_n$ converges strongly to $0$ in $W_0^{s_2,q}(\Omega)$ and consequently in $W_0^{s_1,p}(\Omega)$.
Hence, we infer that \begin{eqnarray*}
\lambda\int_{\Omega}|u_n|^q~dx-
[u_n]^q_{s_2,q}=[u_n]^p_{s_1,p}\rightarrow 0, \text{as $n\rightarrow +\infty$}.
\end{eqnarray*}
Next, using similar argument as the one used in the proof of Proposition \ref{minpositive}, we will reach to a contradiction, which shows that $u_0\not\equiv 0.$ Letting $n\rightarrow \infty$ in relation (\ref{att1}), we deduce that 
\begin{eqnarray*}
[u_0]^p_{s_1,p}+[u_0]^q_{s_2,q}\leq\lambda \int_{\Omega}|u_0|^q~dx.
\end{eqnarray*}
If there is equality in the above relation then $u_0\in\mathcal{N}_{\lambda}$ and $\mathfrak{M}\leq F_{\lambda}(u_0)$. Assume by contradiction that
\begin{equation}\label{att2}
[u]^p_{s_1,p}+[u]^q_{s_2,q}<\lambda \int_{\Omega}|u|^q~dx.
\end{equation} 
Let $t>0$ be such that $tu_0\in\mathcal{N}_{\lambda},$ i.e.,$$t=\Bigg(\frac{\lambda\displaystyle{\int_{\Omega}}|u_0|^q~dx-
[u_0]^q_{s_2,q}}{[u_0]^p_{s_1,p}}\ \Bigg)^{\frac{1}{p-q}}
.$$
We note that $t\in(0,1)$ since $1<t^{p-q}$ (using (\ref{att2})). Finally, since $tu_0\in \mathcal{N}_{\lambda}$ with $t\in(0,1)$ we have 
\begin{eqnarray*}
0<\mathfrak{M}\leq F_{\lambda}(tu_0)&=&\left(\frac{1}{p}-\frac{1}{q}\right)[tu_0]^p_{s_1,p}=t^p\left(\frac{1}{p}-\frac{1}{q}\right)[u_0]^p_{s_1,p}\\
&=&t^p F_{\lambda}(u_0)\\
&\leq & t^p\lim_{n\rightarrow+\infty}\inf F_{\lambda}(u_n)=t^p \mathfrak{M}<\mathfrak{M}~~~~\text{for $t\in(0,1)$.}
\end{eqnarray*}
This is a contradiction, which assures that relation (\ref{att2}) cannot hold and consequently we have $u_0\in \mathcal{N}_{\lambda}$. Hence, $\mathfrak{M}\leq F_{\lambda}(u_0)$ and $ \mathfrak{M}= F_{\lambda}(u_0)$.
\end{proof}

\begin{theorem}\label{constantsign2}
Every $\lambda\in (\lambda_1(s_2,q),+\infty)$ is an eigenvalue  of problem (\ref{Eq1}). Moreover, the associated eigenfunctions $u_{\lambda}$ have a constant sign in $\Omega.$
\end{theorem}
\begin{proof}
Let $u\in \mathcal{N}_{\lambda}$ be such that $F_{\lambda}(u)=\mathfrak{M}$ thanks to Proposition \ref{minatteind}. We show that 
$$\text{$\langle F'_{\lambda}(u), v\rangle=0$\quad  for all \quad $v\in W^{s_2,q}_0(\Omega).$}$$ 
We recall that for $u\in \mathcal{N}_{\lambda}$, we have \begin{eqnarray*}
[u]^p_{s_1,p}+[u]^q_{s_2,q}=\lambda \int_{\Omega}|u|^q~dx.
\end{eqnarray*} 
Let $v\in W^{s_2,q}_0(\Omega).$ For every $\delta$ in some small interval $(-\varepsilon,\varepsilon)$ certainly the function $u+\delta v$ does not vanish identically.~ Let $t(\delta)>0$ be a function such that $t(\delta)(u+\delta v)\in  \mathcal{N}_{\lambda},$ namely 
$$t(\delta)=\Bigg(\frac{\lambda\displaystyle{\int_{\Omega}}|u+\delta v|^q~dx-
[u+\delta v]^q_{s_2,q}}{\displaystyle{[u+\delta v]^p_{s_1,p}}}\ \Bigg)^{\frac{1}{p-q}}
.$$ The function $t(\delta)$ is a composition of differentiable functions, so it is differentiable. The precise expression of $t'(\delta)$ does not matter here. Observe that $t(0)=1.$ The map $\delta\mapsto t(\delta)(u+\delta v) $
defines a curve on $\mathcal{N}_{\lambda}$ along which we evaluate $F_{\lambda}.$ Hence we define ~$\ell: (-\varepsilon,\varepsilon)\rightarrow \mathbb{R}$ as
$$\text{ $\ell(\delta)=F_{\lambda}(t(\delta)(u+\delta v)).$}$$ By construction, $\delta=0$ is a minimum point for $\ell.$ Consequently $$0=\ell'(0)=\langle F'_{\lambda}(t(0)u), t'(0)u+t(0)v\rangle=t'(0)\langle F'_{\lambda}(u),u\rangle+\langle F'_{\lambda}(u),v\rangle=\langle F'_{\lambda}(u),v\rangle.$$ Using the fact that $\langle F'_{\lambda}(u), u\rangle=0$ for $u\in \mathcal{N}_{\lambda}$, we obtained that 
$$\text{$\langle F'_{\lambda}(u),v\rangle=0$\quad  for all \quad $v\in W^{s_2,q}_0(\Omega).$}$$
Next, we show that the eigenfunctions $u_{\lambda}$ has a constant sign in $\Omega.$ If $|u_{\lambda}|\in \cN_{\lambda}$, it follows from Proposition \ref{constantsign} that $u_{\lambda}$ has a constant sign in $\Omega.$ We in fact have that $|u_{\lambda}|\in \cN_{\lambda}$. Suppose by contradiction that $|u_{\lambda}|\notin \cN_{\lambda}$. Then,
\begin{eqnarray*}
[|u_{\lambda}|]^p_{s_1,p}+[|u_{\lambda}|]^q_{s_2,q}<\lambda \int_{\Omega}|u_{\lambda}|^q~dx.
\end{eqnarray*} 
Let $t>0$ be such that $t|u_{\lambda}|\in \cN_{\lambda}$, that is 
$$t=\Bigg(\frac{\lambda\displaystyle{\int_{\Omega}}|u_{\lambda}|^q~dx-
[|u_{\lambda}|]^q_{s_2,q}}{[|u_{\lambda}|]^p_{s_1,p}}\ \Bigg)^{\frac{1}{p-q}}
.$$
We note that $t\in(0,1)$ since $1<t^{p-q}$. Since then $t|u_{\lambda}|\in \mathcal{N}_{\lambda}$ with $t\in(0,1)$, it follows that  
\begin{eqnarray*}
\mathfrak{M}= F_{\lambda}(u_{\lambda})  \le F_{\lambda}(|u_{\lambda}|)\le &=&\left(\frac{1}{p}-\frac{1}{q}\right)[tu_{\lambda}]^p_{s_1,p}=t^p\left(\frac{1}{p}-\frac{1}{q}\right)[u_{\lambda}]^p_{s_1,p}\\\\
&=&t^p F_{\lambda}(u_{\lambda})  =t^p \mathfrak{M}<\mathfrak{M}.
\end{eqnarray*}
This is yields  a contradiction. We conclude that $|u_{\lambda}|\in \cN_{\lambda}$ and the eigenfunctions $u_{\lambda}$ have a constant sign in $\Omega.$ The proof of Theorem \ref{constantsign2} is completed.
\end{proof}

We close this section with  an  analogous result of Theorem \ref{sequence1} for $0<s_1<s_2<1<p<q<\infty.$ it writes as follows:

\begin{theorem}\label{sequence2}
There exists a nondecreasing sequence of critical values $c_k(s_1,s_2,p,q;\rho)$ with associated nonlinear eigenvalues $\lambda_k(s_1,s_2,p,q;\rho)\to +\infty,$ as $k\to +\infty$ and with corresponding eigenfunctions $u_k\in W^{s_2,q}_0(\Omega)$ with $\int_{\Omega}|u_k|^q\ dx=\rho,$ for problem (\ref{Eq1}).
\end{theorem}
\begin{proof}
Let $D_{\rho}(s_1,s_2,p,q)=\{u\in W^{s_2,q}_0(\Omega)~:~\int_{\Omega}|u|^qdx=\rho\}$, and 
$$\Sigma_k(s_1,s_2,p,q)=\{A\subset\Sigma~:~\gamma(A\cap D_{\rho}(s_1,s_2,p,q))\geq k\},$$ where $\Sigma=\{A\subset W^{s_2,q}_0(\Omega):~~A~~\text{closed},~~ A=-A\}.$ Set 
$$
b_k(s_1,s_2,p,q;\rho)=\inf_{A\in\Sigma_k(s_1,s_2,p,q)}\sup_{u\in A}\left(\frac{1}{p}[u]^p_{s_1,p}+\frac 1q[u]^q_{s_2,q}\right)>0.
$$ 
Similar to the proof of Theorem \ref{sequence1}, one shows that:
\begin{enumerate}
    \item the functional $I(u)=\displaystyle{\frac{1}{p}[u]^p_{s_1,p}+ \frac 1q [u]^q_{s_2,q}}$ satisfies the (PS) condition on $D_{\rho}(s_1,s_2,p,q)$, and
    \item if $b=b_k(s_1,s_2,p,q;\rho)=\dots=b_{k+m-1}(s_1,s_2,p,q;\rho),$ then the set $K_b$ of critical points of $I$ at the critical level $b$ has a genus $\gamma(K_b)\geq m.$
\end{enumerate}
\end{proof}

\section{Bifurcation results}\label{S3}

In this section we discuss bifurcation phenomena for problem \eqref{Eq1}. We begin with the following definition.
\begin{definition}\label{bifurcation-theorem}
A real number $\lambda $ is called a bifurcation point of (\ref{Eq1}) if and only if there is a sequence $\{(u_k,\lambda_k)\}_k$ of solutions of (\ref{Eq1}) such that $u_k\not\equiv 0$ and 
$$\lambda_k\rightarrow\lambda, \qquad  \|u_k\|_{s,p}\rightarrow 0,~~~\text{as}~~k\rightarrow\infty,~~\text{if }\quad  0<s_2<s_1<1<q<p<\infty $$
or
$$\lambda_k\rightarrow\lambda, \qquad  \|u_k\|_{s,q}\rightarrow 0,~~~\text{as}~~k\rightarrow\infty,~~\text{if }\quad  0<s_1<s_2<1<p<q<\infty.$$
\end{definition}
We have the following observation. Let consider the following energy functional $G$ defined  on $W^{s_1,p}_0(\Omega)\cap W^{s_2,q}_0(\Omega)$  by:
%$$F(u)=\frac{[u]^p_{s_1,p}}{p}+\frac{[u]^q_{s_2,q}}{q}-\lambda\frac{\|u\|^2_{1,2}}{2},~~\lambda\in\R,$$ where
%$$\|u\|^{2}_{1,2}=\int_{\Omega}|\nabla u|^{2}dx.$$We define the following functional
%$$G(u)=\frac{\frac{[u]^p_{s_1,p}}{p}+\frac{[u]^q_{s_2,q}}{q}}{\frac{\|u\|^q_{q}}{q}}.$$
$$G(u)=\left(\frac{[u]^p_{s_1,p}}{p}+\frac{[u]^q_{s_2,q}}{q}\right)\Big/\frac{\|u\|^q_{q}}{q}.$$
Setting $u=te_1$, where $e_1$ stands for the $L^q$-normalized eigenfunction associated to the first eigenvalue $\lambda_1(s_2,q)$ of the eigenvalue problem \eqref{dirichletlap}. A direct computation shows that
%$$G(te_1)=\frac{\frac{t^{p-q}[u]^p_{s_1,p}}{p}+\frac{[u]^q_{s_2,q}}{q}}{\frac{\|u\|^q_{q}}{q}}.$$
$$G(u)=\left(\frac{t^{p-q}[e_1]^p_{s_1,p}}{p}+\frac{[e_1]^q_{s_2,q}}{q}\right)\Big/\frac{\|e_1\|^q_{q}}{q}=q\left(\frac{t^{p-q}[e_1]^p_{s_1,p}}{p}+\frac{[e_1]^q_{s_2,q}}{q}\right).$$
We distinguish the following two cases:
\begin{enumerate}
    \item Assume that $0<s_1<s_2<1<q<p<\infty.$ We find that $G(te_1)\rightarrow\lambda_1(s_2,q)$ as $t\rightarrow 0,$ which indicates the existence of  bifurcation in $0$ from $\lambda_1(s_2,q).$
    \item Assume that $0<s_2<s_1<1<p<q<\infty.$  We find that $G(te_1)\rightarrow\infty$ as $t\rightarrow 0,$ which indicates there is no bifurcation in $0$ from $\lambda_1(s_2,q).$ We are led to seek for the existence of a bifurcation at infinity, since we have $G(te_1)\rightarrow\lambda_1(s_2,q)$ as $t\rightarrow \infty$.
\end{enumerate}
%\par \bigskip \noindent
The main purpose of this section is to show that any variational $k$-eigenvalue $\lambda_k(s_2,q)$ of problem \eqref{dirichletlap} are bifurcation points for the nonlinear variational eigenvalues $\lambda_k(s_1, s_2, p,q;\rho)$ of problem \eqref{Eq1}. More precisely, we are going to show that
$$ \lambda_k(s_1, s_2, p,q;\rho) \to \lambda_k(s_2,q) \ \quad \text{as}\ \quad \rho \to 0
.$$ 
%We recall that the eigenvalues $\lambda_n(s_2,q)$ are the eigenvalues found in item (4) of Proposition \ref{prop1}.\\
 As in Section \ref{Spec}, we let 
$$D_{\rho}(s_1,s_2,p,q)=\{u\in (W^{s_1,p}_0(\Omega)\cap W^{s_2,q}_0(\Omega))\setminus\{0\}: \int_{\Omega}|u|^q\ dx=\rho\}$$ and
$$\Gamma_{k,\rho}=\{A\subset D_{\rho}(s_1,s_2,p,q): A \ \text{symmetric}, A\ \text{compact}, \gamma(A)\geq k\}.
$$ 
By the definition of $\lambda_k(s_2,q)$ we know that for $\eps>0$ small there is $A_{\eps}\in\Gamma_{k,1}$ such that 
$$
\sup_{\{u\in A_{\eps},~\int_{\Omega}|u|^qdx=1\}}[u]^q_{s_2,q}\leq \lambda_k(s_2,q) +\eps\ .
$$
%(\textbf{ Sur cette page il faut utiliser $n$ comme indice à la place de $k$. })\\
We want to approximate $A_{\eps}$ by a finite-dimensional set. Since $A_{\eps}$ is compact, for every $\delta>0$ there exist a finite number of points $x_1,\dots,x_{n(\delta)}$ such that
\begin{equation}\label{Clarkn1}
    A_{\eps}\subset\bigcup_{i=1}^{n(\delta)}B_{\delta}(x_i).
\end{equation}
Let $E_n=\text{span}\{x_1,\dots,x_{n(\delta)}\},$ and set 
\begin{equation}\label{choose1}
    P_nA_{\eps}:=\{P_nx,~~x\in A_{\eps}\},
\end{equation}
where $P_nx\in E_n$ is such that 
$$\|x-P_nx\|_{s_2,q}=\inf\{\|x-z\|_{s_2,q},~~z\in E_n\}.$$ 
We claim that $\gamma(P_n A_{\eps})\geq k.$ Clearly, $P_n A_{\eps}$ is symmetric and compact. Furthermore, $0\not\in P_n A_{\eps}$. Indeed since $A_{\eps}$ is compact, and $0\not\in A_{\eps}$, there is small ball $B_{\tau}(0)$ such that $A_{\eps}\cap B_{\tau}(0)=\emptyset.$ Now, choose $\delta>0$ in (\ref{Clarkn1}) such that $\delta<\tau/2.$ Then, for $x\in A_{\eps}$ there is $x_i\in E_n$, for some $i \in \{1,\dots,n(\delta)\}$, such that $\|x-x_i\|_{s_2,q}<\delta$, and hence 
$$\|x-P_n x\|_{s_2,q}=\inf\{\|x-z\|_{s_2,q},~~z\in E_n\}\leq \|x-x_i\|_{s_2,q}<\tau/2$$ and thus $P_nA_{\eps}\cap B_{\tau/2}(0)=\emptyset.$\\
Finally, we have to show that $\gamma(P_nA_{\eps})\geq k.$ This is again by approximation: since $\gamma(A_{\eps})\geq k,$ there exists a continuous and odd map $h: A_{\eps}\rightarrow \mathbb{R}^k\setminus\{0\}.$ Then by Tietze extension theorem there exists a continuous and odd map $\tilde{h}: W^{s_2,q}_0(\Omega)\rightarrow \mathbb{R}$ such that $\tilde{h}_{|A_{\eps}}=h.$ By continuity and compactness of $A_{\eps}$ we can conclude that $\tilde{h}_{|P_nA_{\eps}}: W^{s_2,q}_0(\Omega)\rightarrow \mathbb{R}^k\setminus\{0\}.$ Now, again by approximation, we conclude that there is a $n = n(\eps)$ such that 
$$
\sup_{\{u\in P_nA_{\eps}\}}[u]^q_{s_2,q}\leq \lambda_k(s_2,q) + 2\eps\ 
.$$
Finally, note that by homogeneity
$$
\inf_{A\in\Gamma_{k,\rho}}\sup_{u\in A}[u]^q_{s_2,q}=\lambda_k(s_2,q)\,\rho
$$
and hence also
\begin{equation}\label{PnA}
\sup_{\{u\in (\rho\,P_nA_{\eps})\}}[u]^q_{s_2,q}\leq \big(\lambda_k(s_2,q) + 2\eps \big)\,\rho.
\end{equation} 
%\par \bigskip
Recall that by \eqref{ck} we have, for each integer $k > 0$,
$$
c_k(s_1,s_2,p,q;\rho)=\inf_{A\in\Gamma_{k,\rho}}\sup_{u\in A}\Big\{\frac{1}{p}[u]^p_{s_1,p}+\frac 1q[u]^q_{s_2,q}\Big\}
.$$

\subsection{Bifurcation from zero}\label{bafzero}

We assume in this section that  $0<s_2<s_1<1<q<p<+\infty$. We show that  for any $k>0$, problem (\ref{Eq1}) admits a branch of  eigenvalues bifurcating from $(\lambda_k(s_2,q), 0)\in \mathbb{R}^+\times W^{s_1,p}_0(\Omega)$.\\

We prove the following result:

\begin{theorem}\label{bifurcation-from-zero}
For each integer $k>0$ the pair $(\lambda_k(s_2,q), 0)$ is a bifurcation point of problem (\ref{Eq1}).
\end{theorem}

We first prove the following lemma which is the main ingredient for the proof of Theorem \ref{bifurcation-from-zero}.

\begin{lemma}\label{avanbifurc}
 For any integer $k>0$ and  $\rho>0,$ $\eps>0,$ there exists a positive constant $C(\eps)$ such that the following estimate holds: 
$$
|c_k(s_1,s_2,p,q;\rho)- \frac 1q \, \lambda_k(s_2,q)\, \rho|\leq C(\eps)\rho^{p/q} + 2\eps \, \rho
.$$
\end{lemma}

%\todo{{\bf Cette constante ne depend pas seulement de $\eps$}, see my attempt below}\\
%\\
%\textbf{Si elle depend de $\varepsilon$. On utilise simplement l'equivalence des normes dans l'espace de dimension finie $P_kA_\varepsilon$,$~~ k=k(\varepsilon)$.}

\begin{proof}
For any $k>0,$ we  have by the definition of $\lambda_k(s_2,q)$ and $c_k(s_1,s_2,p,q;\rho)$ (see \eqref{min-max} and \eqref{ck}) that $c_k(s_1,s_2,p,q;\rho)\geq \frac 1q\, \lambda_k(s_2,q)\, \rho.$ 
%\par \smallskip \noindent
By \eqref{PnA} we can estimate 
\begin{align*}
    c_k(s_1,s_2,p,q;\rho)&=\inf_{A\in\Gamma_{k,\rho}}\sup_{u\in A}\Big\{\frac{1}{p}[u]^p_{s_1,p}+\frac 1q[u]^q_{s_2,q}\Big\}\\
    &\leq  \sup_{u\in (\rho P_nA_{\eps})}\Big\{\frac{1}{p}[u]^p_{s_1,p}+\frac 1q[u]^q_{s_2,q}\Big\}\\
    &\leq \sup_{u\in (\rho P_nA_{\eps})}\frac{1}{p}[u]^p_{s_1,p}+ \sup_{u\in (\rho\, P_kA_{\eps})}\frac 1q[u]^q_{s_2,q}\\
    &\leq \frac{1}{p}[v]^p_{s_1,p}+\frac 1q(\lambda_n(s_2,q)+2\eps)\rho
\end{align*}
for some $v\in (\rho P_nA_{\eps})$ with $\int_{\Omega}|v|^qdx=\rho.$ Since $ P_nA_{\eps}$ is of finite-dimensional, there exists a positive constant $ C(\eps)$ such that
$$
\big([v]_{s_1,p}\big)^{1/p} \leq C(\eps) \big(\int_{\Omega}|v|^q\ dx\big)^{1/q}.
$$ 
% Indeed, for $v\in (\rho P_kA_{\eps})$, we have $v=\sum_{j=1}^ka_ju_j$, with $a_j\in\R$ (Atomic decomposition) and
% \[
% [v]^p_{s_1,p}\le \sum_{j=1}^ka^p_j[v_j]^p_{s_1,p}\le ?????????\le \lambda_k\sum_{j=1}^ka^{p-q}_j\int_{\Omega}|a_jv_j|^q\ dx.
% \]
% Since $q<p$, it follows using  H\"older inequality that
% \[
% \sum_{j=1}^ka^{p-q}_j\int_{\Omega}|a_jv_j|^q\ dx\le \sum_{j=1}^ka^{p-q}_j\int_{\Omega}|v|^q\ dx\le k^{\frac{p}{q}}\|v\|^{p}_{q}\left(\sum_{j=1}^ka^{p}_j\right)^{\frac{p-q}{p}}%= C(\epsilon,p,q,k)\|v\|^{p}_{q}
% \]
% Hence, it follows that
% $$
% [v]^p_{s_1,p} \leq C(\epsilon,p,q,k) \big(\int_{\Omega}|v|^qdx\big)^{p/q} = C(\epsilon,p,q,k)\, \rho^{p/q}
% .$$
Finally, we get
$$
0\leq c_k(s_1,s_2,p,q;\rho)- \frac 1q\,\lambda_k(s_2,q)\, \rho\leq C(\eps)\rho^{p/q}+2\eps\rho
.$$ 
This completes the proof of Lemma \ref{avanbifurc}.
%that is $$0\leq \lambda^D_k(p,q,\theta)\theta^{-q}-\lambda^D_k(q)\leq C_n(\eps)\theta^{p-q}+\eps. $$
\end{proof}

\begin{proof}[Proof of Theorem \ref{bifurcation-from-zero}]
    We aim at showing that 
    $$\text{$\lambda_k(s_1, s_2, p,q;\rho) \to \lambda_k(s_2,q)$\qquad and\qquad $[u_k]_{s_1,p} \to 0$\quad  as\quad  $\rho\to 0^+.$}$$
 Thanks to Lemma \ref{avanbifurc}, we have  
$$
\frac 1p [u_k]^p_{s_1,p}   
\le C_k(\eps)\rho^{p/q}+2\eps\, \rho 
.$$
Furthermore
$$
\begin{array}{ll}
0 &\le   \lambda_k(s_1,s_2,p,q;\rho)\, \rho - \lambda_k(s_2,q)\rho
\vspace{0.2cm}\\
&=\dis [u_k]^p_{s_1,p} + [u_k]^q_{s_2,q}
- \lambda_k(s_2,q) \rho
\vspace{0.2cm}\\
&
=\dis \frac qp [u_k]^p_{s_1,p} + [u_k]^q_{s_2,q} - \lambda_k(s_2,q) \rho + (1-\frac qp)[u_k]^p_{s_1,p} 
\vspace{0.2cm}\\
&
= \dis q \, c_k(s_1,s_2,p,q;\rho) - \lambda_k(s_2,q)\rho + (1-\frac qp)[u_k]^p_{s_1,p}
\vspace{0.2cm}\\
&
\le  C\, \big(C_k(\eps)\rho^{p/q}+2\eps\, \rho \big).
\end{array}
$$  
Since $\eps > 0$ is arbitrary we get the first claim.
%\par \smallskip

Let us prove that  $[u_k]_{s_1,p}\rightarrow 0$ as $\rho\rightarrow 0^+.$ Letting $v=u_k$ in relation (\ref{Eq3}), we have 
$$
   [u_k]^p_{s_1,p} + [u_k]^q_{s_2,q}=\lambda_k(s_1,s_2,p,q;\rho)\int_{\Omega}|u_k|^{q}~dx.
$$ 
Therefore 
$$ [u_k]^p_{s_1,p} \le \lambda_k(s_1,s_2,p,q;\rho)\int_{\Omega}|u_k|^{q}~dx\\
        \leq C_k\ \rho. $$
Hence $[u_k]_{s_1,p}\rightarrow 0$ as $\rho\rightarrow 0$ and this completes the proof of Theorem \ref{bifurcation-from-zero}.
\end{proof}

\subsection{Bifurcation from infinity}

We assume in this section that  $0<s_1<s_2<1<p<q<+\infty$. The goal in this section is to prove that  there is a branch of eigenvalues bifurcating from $(\lambda_k(s_2,q), +\infty).~~$
For $u\in W^{s_2,q}_0(\Omega),~u\neq 0,$ we set 
$$w=\frac{u}{\|u\|_{s_2,q}^2}.$$  
We have $\|w\|_{s_2,q}=\frac{1}{\|u\|_{s_2,q}}$ and a direct computation shows that for all $x,y\in\Omega$,
\[
\begin{array}{ll}
\|u\|_{s_2,q}^{2(p-q)}| w(x)-w(y)|^{p-2}(w(x)-w(y))&=\frac{1}{\|u\|_{s_2,q}^{2(q-1)}} | u(x)-u(y)|^{p-2}(u(x)-u(y)),
\vspace{5pt}\\
\qquad\qquad| w(x)-w(y)|^{q-2}(w(x)-w(y))&=\frac{1}{\|u\|_{s_2,q}^{2(q-1)}} | u(x)-u(y)|^{q-2}(u(x)-u(y)),\vspace{5pt}\\
\text{ and }\qquad\qquad\qquad\qquad\quad| w(x)|^{q-2}w(x)&=\frac{1}{\|u\|_{s_2,q}^{2(q-1)}} | u(x)|^{q-2}u(x).
\end{array}
\]
Summing up the above  change of variables and using the fact that $u$ is a weak solution of  (\ref{Eq1}), we find that, 
\begin{equation*}\label{em3}
\|u\|_{s_2,q}^{2(p-q)}\cE_{s_1,p}(w,v)+\cE_{s_2,q}(w,v)=\lambda\int_{\Omega}|w|^{q-2}w~v~dx,\qquad\text{ for all $v\in W^{s_2,q}_0(\Omega)$.}
\end{equation*}
This leads  to the following nonlinear eigenvalue problem 
\begin{equation}\label{em5}
\left\{
\begin{array}{rll}
\|w\|_{s_2,q}^{2(q-p)}(-\Delta)^{s_1}_p w+(-\Delta)^{s_2}_q w &=\lambda |w|^{q-2}w~~~~&\text{in\quad $\Omega$} \vspace{0.2cm}\\
w &= \displaystyle 0~~~~~~~~~~~~~~~~&\text{ on \quad $\R^N\setminus\Omega$},
\end{array}
\right.
\end{equation}
where we have used the fact that  $\|u\|^{2(p-q)}_{s_2,q}=\|w\|^{2(q-p)}_{s_2,q}$. \\

We have the following proposition.
\begin{proposition}\label{Prop-Inf1}
 If $(\lambda,0)$ is a bifurcation point of solutions of problem (\ref{em5}) then $\lambda$ is an eigenvalue of the problem \begin{equation}\label{dirichletlap-Infinity}
	\begin{split}
	\quad\left\{\begin{aligned}
(-\Delta)^{s_2}_q u&= \lambda  |u|^{q-2}u && \text{ in \quad $\Omega$}\\
		u &=0     && \text{ on } \quad \R^N\setminus\Omega.
	\end{aligned}\right.
	\end{split}
	\end{equation}
\end{proposition}
\begin{proof}
Since $(\lambda,0)$ is a bifurcation point from zero of solutions of the problem (\ref{em5}), there is a sequence of nontrivial solutions of  problem (\ref{em5}) such that
\[\text{
$\lambda_k\rightarrow\lambda$\qquad and\qquad $\|w_k\|_{s_2,q}\rightarrow 0$\ \ in \  \ \quad as\quad $k\to \infty$}.
\]
We then have by definition of weak solutions that
\begin{equation}\label{pp1}
   \|w_k\|_{s_2,q}^{2(q-p)}\cE_{s_1,p}(w_k,v)+\cE_{s_2,q}(w_k,v)=\lambda_k\int_{\Omega}|w_k|^{q-2}w_kv~dx,\qquad
\end{equation}
\text{ for all $v\in W^{s_2,q}_0(\Omega)$}, ~which equivalent to 
\begin{equation}\label{pp2}
   \|w_k\|_{s_2,q}^{2(p-1)}\cE_{s_1,p}(w_k,v)+\cE_{s_2,q}(u_k,v)=\lambda_k\int_{\Omega}|u_k|^{q-2}u_kv~dx,\qquad
\end{equation}
\text{ for all $v\in W^{s_2,q}_0(\Omega)$}. Since $p>1$, we have  $\|w_k\|_{s_2,q}^{2(p-1)}\cE_{s_1,p}(w_k,v)\to 0$ as $k\to\infty $ for all $v\in W^{s_2,q}_0(\Omega)$. We  pass to the  limit in \eqref{pp2}  using Lemma \ref{positivity} to get 
\begin{equation}
  \cE_{s_2,q}(u,v)=\lambda\int_{\Omega}|u|^{q-2}uv~dx,\quad\text{ for all\quad $v\in W^{s_2,q}_0(\Omega)$}.
\end{equation}
This completes the proof of Proposition \ref{Prop-Inf1}.
\end{proof}
Next, let us consider a small ball $$B_r(0) :=\{~\phi~\in W^{s_2,q}_0(\Omega)\setminus\{0\}/~~~\|\phi\|_{s_2,q}< r~\},$$ and
 the operator $$\cA :=\|\cdot\|_{s_2,q}^{2(q-p)}(-\Delta)^{s_1}_p+(-\Delta)^{s_2}_q : W^{s_2,q}_0(\Omega)\subset W^{s_2,p}_0(\Omega)\longrightarrow W^{-s_2,p'}(\Omega)\subset W^{-s_2,q'}(\Omega).$$
 
\begin{proposition}\label{invert}
There exists $r>0$ small such that the mapping
\[\cA : B_r(0)\subset W^{s_2,q}_0(\Omega)\rightarrow W^{-s_2,q'}(\Omega)\]
is strongly monotone, i.e., there exists $C>0$ such that 
$$\langle \cA(u)-\cA(v), u-v\rangle\geq C\|u-v\|^q_{s_2,q}
, ~~~~\text{ for}~~u,v\in B_r(0),$$ 
with $r>0$ sufficiently small.
\end{proposition}

\begin{proof}
We use the fact that  $(-\Delta)^{s_1}_p$ is strongly monotone on $W^{s_1,p}_0(\Omega)$ (see Lemma \ref{inevec}). By linearity and   H\"older inequality, we have
\begin{eqnarray}\label{bf1}\notag
\langle \cA(u)-\cA(v), u-v\rangle&=&
\| u- v\|^q_{s_2,q}+\left\langle\|u\|_{s_2,q}^{2(q-p)}(-\Delta)^{s_1}_pu-\|v\|_{s_2,q}^{2(q-p)}(-\Delta)_q^{s_2}v, u-v \right\rangle\\ \notag
&=&\| u- v\|^q_{s_2,q}+\|u\|_{s_2,q}^{2(q-p)}\left\langle(-\Delta)^{s_1}_pu-(-\Delta)_p^{s_1}v, u-v\right\rangle\\ \notag
&& \qquad+\left(\|u\|_{s_2,q}^{2(q-p)}-\|v\|_{s_2,q}^{2(q-p)} \right)\left\langle(-\Delta)^{s_1}_pv, u-v\right\rangle\\ \notag
&\geq& \| u- v\|^q_{s_1,q} -\left|\|u\|_{s_2,q}^{2(q-p)}-\|v\|_{s_2,q}^{2(q-p)}\right|\| v\|_{s_2,p}^{p-1}\|u-v\|_{s_2,p}\\
&\geq&\| u- v\|^q_{s_2,q} -\left|\|u\|_{s_2,q}^{2(q-p)}-\|v\|_{s_2,q}^{2(q-p)}\right|C\| v\|_{s_2,q}^{p-1}\| u-v\|_{s_2,q}\\ \notag
&\geq&\| u- v\|^q_{s_2,q}\left(1 -\left|\|u\|_{s_2,q}^{2(q-p)}-\|v\|_{s_2,q}^{2(q-p)}\right|C\| v\|_{s_2,q}^{p-1}\| u-v\|^{1-q}_{s_2,q}\right).
\end{eqnarray}
Moreover, by the Mean Value Theorem,  there exists $\theta\in [0,1]$ such that
\begin{eqnarray*}
\left|\|u\|_{s_2,q}^{2(q-p)}-\|v\|_{s_2,q}^{2(q-p)}\right|&=&\left|\frac{d}{dt}\left(\|u+t(v-u)\|^2_{s_2,q}\right)^{q-p}|_{t=\theta} (v-u)\right|\\
&=&\left|(q-p)\left(\|u+\theta(v-u)\|^2_{s_2,q}\right)^{q-p-1}2\left(u+\theta(v-u),v-u\right)_{s_2,q}\right|\\
&\leq& 2(q-p)\|u+\theta(v-u)\|^{2q-2p-2}_{s_2,q}\|u+\theta(v-u)\|^{q-1}_{s_2,q}\|u-v\|_{s_2,q}\\
&=&2(q-p)\|u+\theta(v-u)\|_{s_2,q}^{2q-2p-1}\|u-v\|_{s_2,q}\\
&\leq & 2(q-p)\|u+\theta(v-u)\|_{s_2,q}^{2q-p}\|u-v\|_{s_2,q}\\ 
&\leq & 2(q-p)\left((1-\theta)\|u\|_{s_2,q}+\theta\|v\|_{s_2,q}\right)^{2q-p}\|u-v\|_{s_2,q}\\
&\leq & 2(q-p) r^{2q-p}\|u-v\|_{s_2,q}.
\end{eqnarray*}
 Substitute the above estimate in (\ref{bf1}), we get
\[
\begin{split}
\langle \cA(u)-\cA(v), u-v\rangle &\geq \|u-v\|^q_{s_2,q}-2(q-p) r^{2q-1}C\|u-v\|^2_{s_2,q}.
\end{split}
\]
Letting $r\to 0,$
the proof of Proposition \ref{invert}  follows.
\end{proof}

We first show the existence of variational eigenvalues of the nonlinear equation (\ref{em5}). 

\begin{theorem}\label{manysolutions}
For a fixed $\rho>0,$ there exists a non-decreasing sequence of eigenvalues $\lambda_k(s_1,s_2,p,q;\rho),$ with corresponding eigenfunctions $w_k\in W^{s_2,q}_0(\Omega) $ for the nonlinear eigenvalue problem (\ref{em5}).
\end{theorem}

We again rely on \cite[Proposition 10.8]{AM} for the proof of Theorem \ref{manysolutions}.
\begin{proof}
 Let define
 $$O_{\rho}(s_1,s_2,p,q)=\{w\in W^{s_2,q}_0(\Omega)~:~\int_{\Omega}|w|^q\ dx=\rho\},$$ and $$\Sigma_{k,\rho}(p,q)=\{A\subset\Sigma~:~\gamma(A\cap O_{\rho}(s_1,s_2,p,q))\geq k\},$$ where $\Sigma=\{A\subset W^{s_2,q}_0(\Omega):~~A~~\text{closed},~~ A=-A\}.$ 
 Set 
 \begin{equation*}\label{dk}
 d_k(s_1,s_2,p,q;\rho)=\inf_{A\in\Sigma_{k,\rho}(p,q)}\sup_{u\in A}\left(\frac{q}{p}\|w\|_{s_2,q}^{2(q-p)}\cE_{s_1,p}(w,w)+\cE_{s_2,q}(w,w)\right)>0.
 \end{equation*} 
We show that:
\begin{enumerate}
    \item the functional $$E(w)=\displaystyle{\frac{q}{p}\|w\|_{s_2,q}^{2(q-p)}\cE_{s_1,p}(w,w)+\cE_{s_2,q}(w,w)}$$ satisfies the (PS) condition on $O_{\rho}(s_1,s_2,p,q)$, and \vspace{0.2cm}
    \item if $d=d_k(s_1,s_2,p,q;\rho)=\dots=d_{k+m-1}(s_1,s_2,p,q;\rho),$ then the set $K_d$ of critical points of $I$ at the critical level $d$ has a genus $\gamma(K_d)\geq m.$
\end{enumerate}

{\bf Proof of 1}.
Let $\{w_j\}\subset O_{\rho}(s_1,s_2,p,q)$ a (PS) sequence, i.e, there is $M>0$ such that
\[\text{
$|E(w_j)|\leq M$\quad and\quad $E'(w_j)\rightarrow 0$ ~~in~~ $W^{-s_2,q'}(\Omega)$~~ as~~ $j \rightarrow \infty.$}
\]
We first show that $\{w_j\}$ is bounded in $O_{\rho}(s_1,s_2,p,q)\subset W^{s_2,p}_0(\Omega).$ Since $w_j\in W^{s_2,q}_0(\Omega),$  it follows using the Poincar\'e inequality and the fact taht $W^{s_2,q}_0(\Omega)\subset W^{s_2,p}_0(\Omega)$, that 
\begin{align*}
    M\geq |E(w_j)|&\geq \frac{q}{p}\|w_j\|_{s_2,q}^{2(q-p)}\cE_{s_1,p}(w_j,w_j)+\frac{1}{C}\int_{\Omega}|w_j|^q\ dx\\
    &\geq \|w_j\|_{s_2,q}^{2q-p}+\frac{\rho}{C}.
\end{align*}
Then $\{w_j\}$ is bounded in $O_{\rho}(s_1,s_2,p,q)\subset W^{s_2,q}_0(\Omega).$ Passing to a subsequence still denoted by $\{w_j\}$, we can assume that there exists $w\in O_{\rho}(s_1,s_2,p,q)\subset W^{s_2,q}_0(\Omega)$ such that $w_j \rightharpoonup w$ in $O_{\rho}(s_1,s_2,p,q)\subset W^{s_2,q}_0(\Omega).$ Now, we show that $w_j$ converges strongly to $w$ in $O_{\rho}(s_1,s_2,p,q)\subset W^{s_2,q}_0(\Omega).$ Since $E'(w_j)\rightarrow 0$ in $W^{-s_2,q'}(\Omega)$ as $j\rightarrow \infty,$ we have 
$$\left\langle E'(w_j)-E'(w), w_j-w\right\rangle\rightarrow 0\quad  \text{ as }\quad j\rightarrow\infty.$$
 %Consider the $C^1$ energy functional associated to equation (\ref{em5}),\\
%$F_{\lambda}: B_r(0)\subset W^{1,q}_0(\Omega)\rightarrow\mathbb{R}$ $$F_{\lambda}(w)=\frac{q}{p}\|w\|_{1,q}^{2(q-p)}\int_{\Omega}|\nabla w|^p~dx+\int_{\Omega}|\nabla w|^q~dx-\lambda\int_{\Omega}|w|^q~dx.$$
 By definition, we have that
 \begin{equation*}
     \begin{split}
         \langle E'(w_j)- &E'(w), w_j-w\rangle
         = p\langle \cA(w_j)-\cA(w), w_j-w_j\rangle.\\
         % &= q\int_{\Omega}\left(\|w_j\|_{1,q}^{2(q-p)}|\nabla w_j|^{p-2}\nabla w_j-\|w\|_{1,q}^{2(q-p)}|\nabla w|^{p-2}\nabla w\right)\cdot \nabla(w_j-w)dx\\
         % &+q\int_{\Omega}\left(|\nabla w_j|^{q-2}\nabla w_j-|\nabla w|^{q-2}\nabla w\right)\cdot \nabla(w_j-w)~dx.
     \end{split}
 \end{equation*}
Thanks to Proposition \ref{invert}, it follows that
 \begin{equation*}
     \begin{split}
         \langle E'(w_j)- E'(w), w_j-w\rangle &\geq C \|w_j-w\|^q_{s_2,q}.
     \end{split}
 \end{equation*}
 Therefore $\|w_j-w\|_{s_2,q}\rightarrow 0$ as $j\rightarrow +\infty$ and $w_j$ converges strongly to $w$ in $W^{s_2,q}_0(\Omega).$
% \par \medskip
 
  The proof of $2.$  is similar to the last part of the proof of Theorem \ref{sequence1}.
 %\todo{To check carefully till the end}
\end{proof}
%\par \medskip
\begin{theorem}\label{bifurcation-from-infinity}
For each integer $k>0$ the pair $(\lambda_k(s_2,q,\rho), +\infty)$ is a bifurcation point of problem (\ref{Eq1}). 
\end{theorem}
The proof of Theorem \ref{bifurcation-from-infinity} will follow immediately from the following remark, and the proof that $(\lambda_k(s_2,q,\rho),0)$ is a bifurcation point of (\ref{em5}), which will be shown in Theorem \ref{binfnon} below.

\begin{remark}\label{RMK}
With the change of variable $\frac{u}{\|u\|^2_{s_2,q}}$, we have that the pair $(\lambda_k(s_2,q,\rho),+\infty)$ is a bifurcation point for the problem (\ref{Eq1}) if and only if the pair $(\lambda_k(s_2,q,\rho),0)$ is a bifurcation point for the problem (\ref{em5}).
\end{remark}
%\par \medskip
Before we proceed to the proof of Theorem \ref{binfnon} below, we show the following lemma.
\begin{lemma}\label{bifurcainfinity}
Let $0<s_1<s_2<1<p<q<+\infty$. For any integer $k>0$ and  $\rho>0,$ $\eps>0,$ there exists a positive constant $D(\eps)$ such that the following estimate holds: $$|d_k(s_1,s_2, p,q;\rho)-\lambda_k(s_2,q,\rho)|\leq (D(\eps)+\eps)\rho^{\frac{2q-p}{p}}
$$ 
where $d_k(s_1,s_2,p,q;\rho)$ is given by \eqref{dk}, and 
$$\dis \lambda_k(s_2,q,\rho)=\inf_{A\in\Gamma_{k,\rho}}\sup_{u\in A}\cE_{s_2,q}(w,w)=\lambda_k(s_2,q,\rho)\rho.
$$
\end{lemma}

\begin{proof}
For any $k>0,$ we clearly have $d_k(s_1,s_2,p,q;\rho)\geq \lambda_k(s_2,p,\rho).$ As in (\ref{choose1}), we choose $P_nA_{\eps}$ such that 
$$\sup_{\{w\in P_nA_{\eps},~\int_{\Omega}|w|^qdx=1\}}\cE_{s_2,q}(w,w)\leq \lambda_k(s_2,q,\rho)+\eps
$$ and so 
$$
\sup_{\{w\in P_nA_{\eps,\rho},~\int_{\Omega}|w|^qdx=\rho\}}\cE_{s_2,q}(w,w)\leq (\lambda_k(s_2,q,\rho)+\eps)\rho,
$$ 
where $P_nA_{\eps,\rho}=\{w\in P_nA_{\eps}:~~\int_{\Omega}|w|^qdx=\rho\}.$
Then 
\begin{align*}
    d_k(s_1,s_2,p,q;\rho)&=\inf_{A\in\Gamma_{k,\rho}}\sup_{w\in A}\Big\{\frac{q}{p}\|w\|_{s_2,q}^{2(q-p)}\cE_{s_1,p}(w,w)+\cE_{s_2,q}(w,w)\Big\}\\
    &\leq  \sup_{w\in P_nA_{\eps,\rho}}\Big\{\frac{q}{p}\|w\|_{s_2,q}^{2(q-p)}\cE_{s_1,p}(w,w)+\cE_{s_2,q}(w,w)\Big\}\\
    &\leq \sup_{w\in P_nA_{\eps,\rho}}\frac{q}{p}\|w\|_{s_2,q}^{2(q-p)}\cE_{s_1,p}(w,w)+ \sup_{w\in P_nA_{\eps,\rho}}\cE_{s_2,q}(w,w)\\
    &\leq\frac{q}{p}\|v\|_{s_2,q}^{2(q-p)}\cE_{s_2,q}(v,v)+(\lambda_k(s_2,q)+\eps)\rho~~\text{since}~~p<q,\\
    &\leq\frac{q}{p}\|v\|_{s_2,q}^{2q-p}+(\lambda_k(s_2,q)+\eps)\rho
\end{align*}
for some $v\in P_nA_{\eps,\rho}$ with $\int_{\Omega}|v|^q\ dx=\rho.$ Since $ P_nA_{\eps}$ is finite-dimensional, there exists a positive constant $D_k(\eps)$ such that
$\cE_{s_2,q}(v,v)\leq D_k(\eps) (\int_{\Omega}|v|^qdx)^{p/q}=D_k(\eps)\rho^{q/p}$ and 
$$\|v\|_{s_2,q}^{2q-p}\leq D_k(\eps)\rho^{\frac{2q-p}{p}}.$$
Finally, we get
$$0\leq d_k(s_1,s_2,p,q;\rho)-\lambda_k(s_2,q,\rho)\leq D_k(\eps)\rho^{\frac{2q-p}{p}}+\eps\rho\leq (D_k(\eps)+\eps)\rho^{\frac{2q-p}{p}}$$ since $\frac{2q-p}{p}>1.$
\end{proof}
\begin{remark}\label{bifurcation-infty}
    We recall that the $k$-th eigenvalue of equation (\ref{em5}) satisfies 
    $$
    \tilde{\lambda}_k(s_1,s_2,p,q;\rho)\rho =\|w\|^{2(q-p)}_{s_2,p}\cE_{s_1,p}(w,w)+\cE_{s_2,q}(w,w),~~\text{\ with\ }~~\rho=\int_{\Omega}|w|^q\ dx.
    $$ 
    So, proceeding as in Theorem \ref{bifurcation-from-zero} one obtains that $\tilde \lambda_k(s_1,s_2,p,q;\rho)\to {\lambda}_k(s_2,q)$ as $\rho\to 0^+.$
\end{remark}
%\par \medskip
\begin{theorem}\label{binfnon}
For any   $k>0$, the pair $(\lambda_k(s_2,q),0)$ is a bifurcation point of problem (\ref{em5}).
\end{theorem}
\begin{proof}
In order to prove Theorem \ref{binfnon}, it suffices to prove that $\tilde{
\lambda}_k(s_1,s_2,p,q;\rho)\rightarrow\lambda_k(s_2,q,\rho)$ and $\|w_k\|_{s_2,q}\rightarrow 0$ as $\rho\rightarrow 0^+.$ 
%where $\tilde{\lambda}^D_k(p,q)$ denotes the $k$-th eigenvalue of problem (\ref{em5}) found in Theorem \ref{manysolutions}.
The fact that $\tilde{
\lambda}_k(s_1,s_2,p,q;\rho)\rightarrow\lambda_k(s_2,q)$ as $\rho\rightarrow 0^+$ follows from Lemma \ref{bifurcainfinity} and Remark \ref{bifurcation-infty}. 

It remains to prove that~~ $\|w_k\|_{s_2,q}\rightarrow 0$ as $\rho\rightarrow 0^+.$ For any $k>0,$ we have
\begin{equation*}
    \begin{split}
       \|w_k\|_{s_2,q}^{2(q-p)}\cE_{s_1,p}(w_k,w_k) + \cE_{s_2,q}(w_k,w_k) &=\tilde{\lambda}_k(s_1,s_2,p,q;\rho)\int_{\Omega}|w_k|^q~dx
       \\
       & \le C_k \int_{\Omega}|w_k|^q~dx= C_k\, \rho \to 0 \ , \ \hbox{ as } \ \rho \to 0. 
       \end{split}
\end{equation*}
Therefore $\|w_k\|_{s_2,q} \to 0$, and since $p<q,$ by the H\"older inequality there exists a positive constant $C_1$ such that $\cE_{s_1,p}(w_k,w_k)\leq C_1\|w_k\|^p_{s_2,q}$, and so also $\|w_k\|_{s_2,p} \to 0$. This completes the proof.
\end{proof}
%\par \bigskip

\section{Multiplicity results: Proof of Theorem \ref{Mainresult2} }\label{Multi}

In this section, we prove a multiplicity result. We show that for fixed $\lambda \in (\lambda_k(s_2,q), \lambda_{k+1}(s_2,q))$ there exist at least $k$ pairs of eigenfunctions \
$\pm u_k^{\lambda}, k = 1,\dots,k, $ such that $(\lambda,\pm u_k^\lambda)$ solves problem \eqref{e2}, i.e. 
$$
\lambda = \lambda_1(s_1,s_2,p,q,\rho_1) = \dots = \lambda_k(s_1,s_2,p,q;\rho_k)\ , \ \hbox{ with } \quad  \int_\Omega |u_k^\lambda|^q\ dx =\rho_k. 
$$

For convenience of the readers, we recall the statement of Theorem \ref{Mainresult2} here.
\begin{theorem}\label{t3}
Let $0<s_2<s_1<1<q<p<\infty$ or $0<s_1<s_2<1<p<q<\infty$, and suppose that $\lambda\in(\lambda_k(s_2,q),\lambda_{k+1}(s_2,q))$ for any $k>0$ and $\rho>0.$ Then equation (\ref{Eq1}) has at least $k$ pairs of nontrivial solutions.
\end{theorem}

\noindent We distinguish  two cases whether: $0<s_2<s_1<1<q<p<\infty$ or $0<s_1<s_2<1<p<q<\infty$.

\subsection{The case:  $0<s_1<s_2<1<p<q<\infty$}

\begin{proof}[{\bf Proof of Theorem \ref{t3}}] The proof in this case  relies  on variational methods and  we will make use of \cite[Proposition 10.8]{AM}. We consider the functional $J_{\lambda}: W^{s_2,q}_0(\Omega)\backslash\{0\}\rightarrow\mathbb{R}$ associated to the problem (\ref{Eq1}) defined by
$$J_{\lambda}(u)=\frac{q}{p}[u]^p_{s_1,p}+[u]^q_{s_2,q}-\lambda\int_{\Omega}|u|^q~dx.$$
Note that  the functional $J_{\lambda}$ is not bounded from below on $W^{s_2,q}_0(\Omega)$. So, we consider  the  constraint set $\mathcal{N}_{\lambda}$  on which we minimize the functional $J_{\lambda}.$ We recall that the constraint set $\mathcal{N}_{\lambda}$ is defined by
$$\mathcal{N}_{\lambda}:=\{u\in W^{s_2,q}_0(\Omega)\backslash\{0\}:~\langle J'_{\lambda}(u),u\rangle=0\}.$$
We have  that $J_{\lambda}(u)=(\frac{1}{p}-\frac{1}{q})[u]^p_{s_1,p}>0$ on $\mathcal{N}_{\lambda}$and it is even and bounded from below on $\mathcal{N}_{\lambda}.$
We show that every Palais-Smale (PS) sequence for $J_{\lambda}$ has a converging subsequence on $\mathcal{N}_{\lambda}.$ 

Let $\{u_n\}_{n}$ be a (PS) sequence, i.e, $|J_{\lambda}(u_n)|\leq C$, for all $n$, for some $C>0$ and $J'_{\lambda}(u_n)\rightarrow 0$ in $W^{-s_2,q'}(\Omega)$ as $n\rightarrow +\infty,$ with $\frac{1}{q}+\frac{1}{q'}=1.$  
We first show that the sequence $(u_n)_{n\geq 0}$ is bounded on $\mathcal{N}_{\lambda}$.
Suppose that $(u_n)_{n\geq 0}$ is not bounded, so $[u_n]^q_{s_2,q}\rightarrow +\infty$ as $n\rightarrow +\infty.$ Since $J_{\lambda}(u_n)=(\frac{1}{p}-\frac{1}{q})[u_n]^p_{s_1,p},$ we have $[u_n]^p_{s_1,p}\leq C.$ On $\mathcal{N}_{\lambda}$, we have
 \begin{equation}\label{ml1}
 0<[u_n]^p_{s_1,p}=\lambda\int_{\Omega}|u_n|^q~dx-[u_n]^q_{s_2,q}.
 \end{equation}
  Set $v_n=\frac{u_n}{\|u_n\|_q}$ then $[v_n]^q_{s_2,q}< \lambda$ (using (\ref{ml1}) and hence $v_n$ is bounded in $W^{s_2,q}_0(\Omega).$ Therefore, there exists $v_0\in W^{1,q}_0(\Omega)$ such that $v_n\rightharpoonup v_0$ in $W^{s_2,q}_0(\Omega)$ and $v_n\rightarrow v_0$ in $L^q(\Omega)$ and $v_n\rightarrow v_0$ in $L^p(\Omega)$ (since $p<q$). Dividing (\ref{ml1}) by $\|u_n\|^p_q,$ we have
 $$\frac{\lambda\displaystyle{\int_{\Omega}}|u_n|^q~dx-[u_n]^q_{s_2,q}}{\|u_n\|^p_q}=[v_n]^p_{s_1,p}\rightarrow 0, \quad \text{ as }\quad n\to \infty.$$ 
 This because,  $\lambda\displaystyle{\int_{\Omega}}|u_n|^q~dx-[u_n]^q_{s_2,q}\le(\frac{1}{p}-\frac{1}{q})^{-1}|J_{\lambda}(u_n)|\leq C$ and $\|u_n\|^p_q\rightarrow +\infty.$ Now, since $v_n\rightharpoonup v_0$ in $W^{s_2,q}_0(\Omega)\subset W^{s_1,p}_0(\Omega),$ we infer that 
 $$\int_{\Omega}|v_n|^p\ dx\le C[v_0]^p_{s_1,p}\leq C\liminf_{n\rightarrow +\infty}[v_n]^p_{s_1,p}=0.$$ 
 Consequently, $v_0=0.$  This is a contradiction, since $\|v_n\|_q=1.$ Thus $\{u_n\}_{n}$ is bounded on $\mathcal{N}_{\lambda}.$ 
 
 Next, we show that $u_n$ converges strongly to $u$ in $W^{s_2,q}_0(\Omega).$ We have 
 $$\text{$\displaystyle{\int_{\Omega}}|u_n|^{q-2}u_n~dx\rightarrow\displaystyle{\int_{\Omega}}|u|^{q-2}u~dx$\quad  as \quad $n\rightarrow+\infty$.}$$ Since $J_{\lambda}'(u_n)\rightarrow 0$ in $W^{-1,q'}(\Omega)$ and  $u_n\rightharpoonup u$ in $W^{1,q}_0(\Omega),$ it follows that $J_{\lambda}'(u_n)(u_n-u)\rightarrow 0$ and $J_{\lambda}'(u)(u_n-u)\rightarrow 0$ as $n\rightarrow +\infty.$ A straightforward computations shows that for $1<p<+\infty$, 
\begin{align*}
\langle(-\Delta)^{s_1}_pu_n-(-\Delta)^{s_1}_pu, u_n-u\rangle
    \geq 0.
\end{align*}
 Therefore, 
 \begin{eqnarray*}
 \langle J'_{\lambda}(u_n)- J'_{\lambda}(u), u_n-u\rangle &=& q\left[\langle(-\Delta)^{s_1}_pu_n-(-\Delta)^{s_1}_pu, u_n-u\rangle\right]\\
 &+& q\left[\langle(-\Delta)^{s_2}_qu_n-(-\Delta)^{s_2}_qu, u_n-u\rangle\right]\\ &-& \lambda q\left[\int_{\Omega}\left(| u_n|^{q-2}u_n-| u|^{q-2} u\right)\cdot (u_n-u)~dx\right]\\
 &\geq & q\left[\langle(-\Delta)^{s_2}_qu_n-(-\Delta)^{s_2}_qu, u_n-u\rangle\right]\\
 &-& \lambda q\left[\int_{\Omega}\left(| u_n|^{q-2}u_n-| u|^{q-2} u\right)\cdot (u_n-u)~dx\right]. 
 \end{eqnarray*}
 Using Lemma \ref{inevec}, it follows that
 \begin{equation*}
     \begin{split}
         \langle J'_{\lambda}(u_n)- J'_{\lambda}(u), u_n-u\rangle &\geq C\|u_n-u\|^q_{s_2,q}-\lambda q\left[\int_{\Omega}\left(| u_n|^{q-2}u_n-| u|^{q-2} u\right)\cdot (u_n-u)~dx\right].
     \end{split}
 \end{equation*}
 Therefore $\|u_n-u\|_{s_2,q}\rightarrow 0$ as $n\rightarrow +\infty$ and $u_n$ converges strongly to $u$ in $W^{s_2,q}_0(\Omega).$
 \\
 \\
 Let $\Sigma=\{A\subset\mathcal{N}_{\lambda}:~A~\text{closed}~\text{and}~-A=A\}$ and $\Gamma_j=\{A\in\Sigma:~\gamma(A)\geq j\}.$\\  We show that $\Gamma_j\neq\emptyset,$ for $j\in \{1,\dots,k\}$.

 Let $\lambda\in (\lambda_j(s_2,q),\lambda_{j+1}(s_2,q))$ and choose  $S^{\eps}_j\in\Sigma\cap\{\int_{\Omega}|u|^q~dx=1\}$   such that $$\sup_{v\in S^{\eps}_j}[v]^q_{s_2,q}\leq \lambda_j(s_2,q)+\eps,~~\eps:=\frac{\lambda-\lambda_j(s_2,q)}{2}.$$
 Then, for $v\in S^{\eps}_j$
 we set
$$\rho(v)=\left[\frac{[v]^p_{s_1,p}}{\lambda\int_{\Omega}|v|^q~dx-[v]^q_{s_2,q}}\right]^{\frac{1}{q-p}},$$ with 
\begin{align*}
    \lambda\int_{\Omega}|v|^q~dx-[v]^q_{s_2,q}&\geq \lambda\int_{\Omega}|v|^q~dx-(\lambda_j(s_2,q)+\eps)\int_{\Omega}| v|^q~dx\\
    &=(\lambda-\lambda_j(s_2,q)-\eps)\int_{\Omega}| v|^q~dx\\
    &=[\lambda-\lambda_j(s_2,q)-(\frac{\lambda-\lambda_j(s_2,q)}{2})]\int_{\Omega}| v|^q~dx\\
    &=\frac{\lambda-\lambda_j(s_2,q)}{2}\int_{\Omega}| v|^q~dx>0,~~\text{for all}~~v\in S^{\eps}_j.
\end{align*} Hence, $\rho(v)v\in \mathcal{N}_{\lambda},$ and then $\rho(S^{\eps}_j)\in\Sigma,$ and $\gamma(\rho(S^{\eps}_j))=\gamma(S^{\eps}_j)=j$ for $1\leq j\leq k.$

It is then standard \cite[Proposition 10.8]{AM} to conclude that $$\sigma_{\lambda,j}=\inf_{A \in \Gamma_j}\sup_{u\in A} J_{\lambda}(u),~~1\leq j\leq k,~~\text{for any }~k\in\mathbb{N}^*$$ yields $k$ pairs of nontrivial critical points for $J_{\lambda},$ which gives rise to $k$ nontrivial solutions of problem (\ref{Eq1}).
\end{proof}

\subsection{The case:  $0<s_2<s_1<1<q<p<\infty$}

\begin{proof}[{\bf Proof of Theorem \ref{t3}}] In this case, we will rely on the following theorem due to Clark \cite{Cl}.
%\textbf{Theorem} \textup{(Clark, \cite{Cl}) }.\\

\begin{theorem}\label{cl1}
Let $X$ be a Banach space and $F\in C^1(X,\mathbb{R})$ satisfying the (PS)-condition with $F(0)=0.$ Let $$\text{$\Gamma_k =\{~A\in\Sigma~:~\gamma(A)\geq k~\}$\quad with\quad $\Sigma= \{~A\subset X~;~A=-A~\text{and}~A~\text{closed}~ \}$.}$$ If 
$c_k=\inf\limits_{A\in \Gamma_k}\sup\limits_{u\in A}F(u)\in (-\infty, 0),$
then $c_k$ is a critical value.
\end{theorem}

We consider the $C^1$ functional $J_{\lambda}: W^{s_1,p}_0(\Omega)\subset W^{s_2,q}_0(\Omega)\rightarrow\mathbb{R},$ $$J_{\lambda}(u)=\frac{q}{p}[u]^p_{s_1,p}+[u]^q_{s_2,q}-\lambda\int_{\Omega}|u|^q~dx.$$ 

Let $\Gamma_k=\{A\subset W^{s_2,q}_0(\Omega)\backslash\{0\},~~A~~\text{compact},~~A=-A,~~\gamma(A)\geq k\},$ and for $\eps>0$ small let $A_{\eps}\in\Gamma_k$ such that $$\sup_{\{u\in A_{\eps},~\int_{\Omega}|u|^qdx=1\}}[u]^q_{s_2,q}\leq \lambda_k(s_2,q)+\eps.$$
We would like to show that 
\begin{equation}\label{ml2}
-\infty<\alpha_{\lambda,k}=\inf_{A\in\Gamma_k}\sup_{u\in A}J_{\lambda}(u)
\end{equation} 
are critical values for $J_{\lambda}.$ We clearly have that $J_{\lambda}(u)$ is  an even functional for all $u\in W^{s_1,p}_0(\Omega)$, and also $J_{\lambda}$ is bounded from below on $W^{s_1,p}_0(\Omega)$ since $J_{\lambda}$ is coercive on $W^{s_1,p}_0(\Omega)$. 

We show that $J_{\lambda}(u)$ satisfies the (PS) condition. Let $\{u_n\}$ be a Palais-Smale sequence, i.e., $|J_{\lambda}(u_n)|\leq M$ for all $n,$ $M>0$ and $J_{\lambda}'(u_n)\rightarrow 0$ in $W^{-s_1,p'}(\Omega)$ as $n\rightarrow+\infty.$  We first show that $\{u_n\}$ is bounded in $W^{1,p}_0(\Omega).$ We have 
\begin{eqnarray*}
M&\geq & |C\|u_n\|_{s_1,p}^p-C'\|u_n\|^q_{s_1,p}|\geq | C\|u_n\|_{s_1,p}^{p-q}-C'| \|u_n\|_{s_1,p}^q,
\end{eqnarray*}
and so $\{u_n\}$ is bounded in $W^{s_1,p}_0(\Omega).$
Therefore, $u\in W^{s_1,p}_0(\Omega)$ exists such that, up to subsequences that we will denote by $(u_n)_n$ we have $u_n\rightharpoonup u$ in $W^{s_1,p}_0(\Omega)$ and $u_n\rightarrow u$ in $L^q(\Omega).$ Arguing as in Part 1, we obtain that $\|u_n-u\|_{s_1,p}\rightarrow 0$ as $n\rightarrow +\infty$, and so $u_n$ converges to $u$ in $W^{s_1,p}_0(\Omega)\subset W^{s_2,q}_0(\Omega) .$

%\par \medskip
As in section \ref{S3}, we approximate $A_{\eps}$ by a finite-dimensional set. Next, we show that there exists sets $D^{\eps}$ of genus greater or equal to $k$ such that $\sup\limits_{u\in D^{\eps}}J_{\lambda}(u)<0.$
For any $t\in (0,1)$, we define the set $D^{\eps}(t):=t\cdot (P_nA_{\eps})$ and so $\gamma(D^{\eps}(t))=\gamma(P_nA_{\eps})\geq k .$ We have, for any $t\in (0,1)$
\begin{eqnarray}\notag
\sup\limits_{u\in D^{\eps}}J_{\lambda}(u) &= & \sup\limits_{u\in P_nA_{\eps}}J_{\lambda}(tu)\\\notag
&\leq& \sup\limits_{u\in P_nA_{\eps}}\left\{\frac{qt^p}{p}[u]^p_{s_1,p}+t^q[u]^q_{s_2,q}-\lambda t^q\int_{\Omega}|u|^qdx\right\}\\
&\leq & \sup\limits_{u\in P_nA_{\eps}}\left\{\frac{qt^p}{p}c(n)^p\|u\|^p_{s_2,q}+t^q(\lambda_k(s_2,q)+\eps-\lambda)\right\}<0\notag
\end{eqnarray}
for $t>0$ sufficiently small.
Finally, we conclude that $\alpha_{\lambda,k}$ are critical values for $J_{\lambda}$ thanks to Clark's Theorem.
 \end{proof}
 
\addcontentsline{toc}{section}{References}

\bibliography{refs}

\vspace{1.5cm}
\noindent
Emmanuel Wend-Benedo Zongo\\
Université Toulouse III - Paul Sabatier,  Institut de mathématiques de Toulouse, France \\
\noindent
\textit{Email address}: ~\texttt{emmanuel.zongo@math.univ-toulouse.fr}

\vspace{1cm}
\noindent
Pierre Aime Feulefack\\
{University of Pennsylvania (Upenn). Philadelphia, PA 19104, United States}\\
\noindent
\textit{Email address}: ~\texttt{pierre.feulefack@aims-cammeroon.org}
\end{document}